\documentclass[12pt]{article}

\usepackage[english]{babel}
\usepackage{indentfirst}
\usepackage{latexsym}
\usepackage[colorlinks=true, bookmarksnumbered=true, bookmarksopen=true,
bookmarksopenlevel=3, pdfstartview=FitH, linkcolor=cyan, pdfmenubar=true,
pdftoolbar=true, bookmarks=true,citecolor=cyan, urlcolor=magenta,
filecolor=cyan,menucolor=black,plainpages=false,pdfpagelabels]{hyperref}
\usepackage[lmargin=2cm,tmargin=2cm,bmargin=2cm,rmargin=2cm]{geometry}
\usepackage{amsbsy, amsmath,amsfonts,amssymb,amsthm}
\usepackage{times}
\usepackage{graphicx}
\usepackage{amsmath}
\usepackage{amsfonts}
\usepackage{amssymb}
\usepackage[english]{babel}
\usepackage[latin1]{inputenc}
\usepackage[T1]{fontenc}
\usepackage{amsmath,amssymb,graphics,amsthm}
\usepackage[usenames]{color}

\numberwithin{equation}{section}
\newtheorem{prop}{Proposition}[section]
\newtheorem{defi}[prop]{Definition}
\newtheorem{teo}[prop]{Theorem}

\newtheorem{obs}[prop]{Remark}
\newtheorem{lema}[prop]{Lemma}

\newcommand\supp{\mathop{\rm supp}}
\newcommand\reallywidehat[1]{\savestack{\tmpbox}{\stretchto{  \scaleto{    \scalerel*[\widthof{\ensuremath{#1}}]{\kern-.6pt\bigwedge\kern-.6pt}    {\rule[-\textheight/2]{1ex}{\textheight}}  }{\textheight}}{0.5ex}}\stackon[1pt]{#1}{\tmpbox}}

\begin{document}

\title{Besov mixed-Morrey spaces:
On an Application 
\\
to the Navier-Stokes Equations}
\author{
{ Leithold L. Aurazo-Alvarez$^{1}$}{\thanks{
{E-mail address: aurazoall@gmail.com}
\newline
{L.L. Aurazo-Alvarez was supported by FAPERJ (Grants 200.140/200.141/2024) and Federal University of Rio de Janeiro, RJ, Brazil.
 }}} , 
%EndAName
\, {Wladimir Neves $^{2}$}{\thanks{{Corresponding author. }
\newline
{E-mail address: wladimir@im.ufrj.br}\newline
%,lcff@ime.unicamp.br (L.C.F. Ferreira).}\newline
{Wladimir Neves has received research grants from FAPERJ (Cientista do Nosso Estado) through the grant E-26/204.171/2024, 
and also from CNPq through the grant  313005/2023-0, 406460/2023-0.
 %L.C.F. Ferreira was supported by FAPESP and CNPq, Brazil.
 }}} 
 \\[5pt]
%EndAName
{\small $^{1,2}$ Instituto de Matem\'atica, Universidade Federal do Rio de Janeiro } 
\\[5pt]
{\small Cidade Universit\'aria, CEP 21945-970, Rio de Janeiro, Brazil.}
}
\date{\today}
\maketitle
%%%%%%%%%%%%%%%
%%%%%%%%%%

\begin{abstract}
In this paper we introduce two new classes of functional spaces, namely, 
Besov mixed-Morrey spaces and Fourier-Besov mixed-Morrey spaces, and 
then we establish some basic properties for these classes.
Moreover, we explore the $d$-dimensional incompressible Navier-Stokes 
equations in this context, by mean of a Bony's paraproduct approach, in order 
to get the global well-posedness for small initial data. These results provides 
a new class of initial data with a sort of anisotropy in relation to its spatial variables or frequency variables.
 
 \medskip
 

{\small \medskip\noindent\textbf{Key:} 
Well-posedness; Navier-Stokes equations; Critical spaces; Besov mixed-Morrey spaces; Fourier-Besov mixed-Morrey spaces; Mixed-Morrey spaces; Large initial data.}
\end{abstract}

\renewcommand{\abstractname}{Abstract}

\tableofcontents

%%%%%%%%%%%%%
\section{Introduction}
%%%%%%%%%%%%%

In this work we propose two classes of functional spaces, 
namely, Besov mixed-Morrey spaces and Fourier-Besov mixed-Morrey spaces,
containing functions with a sort of anisotropy, 
in order to explore the Navier-Stokes equations. For instance, for homogeneous Besov mixed-Morrey spaces, 
it is possible to consider the initial velocity as in Remark \ref{initial-data-example}. 
We build upon key ideas and computations from the seminal work of Kozono and 
Yamazaki \cite{Kozono-Yamazaki-semilinear-1994}, which addresses homogeneous 
Besov-Morrey spaces, as well as from the foundational contribution by Almeida, 
Ferreira, and Lima \cite{Almeida-Ferreira-Lima-uniform-2017}, 
concerning homogeneous Fourier-Besov-Morrey spaces.
To the best of the authors' knowledge, this is the first work to explore these 
mixed spaces and their application to the Navier-Stokes equations.

\medskip
Besov spaces of functions have played a relevant role in the study of various kinds 
of partial differential equations. In particular, they provide a finer functional framework 
than Sobolev spaces for establishing the existence of local or global solutions to initial 
value problems with low-regularity data. For instance, in the case of the $d$-dimensional 
incompressible Navier-Stokes equations \eqref{BQ}, hereafter referred to as the "(N-S) 
Equations", they are useful for analyzing the regularity, spatial decay, 
or asymptotic behavior of the velocity field over time.
Moreover, the new class of function spaces considered in this paper is more suitable 
for capturing local or directional regularity, with potential applications to flows that are 
slow compared to the speed of sound, as encountered in oceanography and meteorology. 
It is worth recalling that real fluid flows often exhibit rough structures or intermittent behavior.

%%%%%%%%%%%%%%%
\subsection{Contextualization}
%%%%%%%%%%%%%%%

In last years some researches were involved to study the Navier-Stokes equations by considering mixed spaces. 
In particular, mixed-Lebesgue spaces, Besov spaces based on mixed-Lebesgue spaces, Fourier-Besov spaces based on mixed-Lebesgue spaces or Besov spaces based on mixed-Lorentz spaces.
In the work \cite{Phan-mixed-2020} T. Phan proved the global existence of a unique mild solution, for small 
initial data $u_{0}\in L^{\vec{p}}(\mathbb{R}^{d})$, with $\vec{p}= (p_1, \ldots, p_d) \in [2,\infty)^{d}$, extending the previous result by Kato \cite{Kato-strong-1984} 
(for mixed-$L^{\vec{p}}(\mathbb{R}^{d})$ spaces, see \cite{Benek-Panzone-mixed-Lp-1961}). The novelty for this approach was the 
boundedness for the Riesz transform on $L^{\vec{p}}(\mathbb{R}^{d})$ and thus for the Leray operator. The solution spaces allows 
functions which decay to zero as $\mid x\mid$ tends to infinity, by considering different rates according to each spatial direction. 
Moreover, this result provides the persistence of anisotropic properties along the evolution of the system.  Khai and Tri \cite{Khai-Tri-Sobolev-Lorentz-Navier-Stokes-2014} proved the existence of a unique global mild solution, for small initial data in the Besov spaces based on mixed-Lorentz spaces $L^{\vec{q},\vec{r}}(\mathbb{R}^{d})$, or more compactly mixed-norm Besov-Lorentz spaces, that is, $u_{0}\in {\cal \dot{B}}^{-\frac{2}{p}}_{L^{\vec{q},\vec{r}},p}(\mathbb{R}^{d})$, for $\vec{r}\in [1,\infty)^{d}$, $\vec{q}\in (2,\infty)^{d}$ and $2<p<\infty$ such that $\frac{-2}{p}=-1+\displaystyle{\sum_{i=1}^{d}}\frac{1}{q_{i}}$, which implies $\displaystyle{\sum_{i=1}^{d}}\frac{1}{q_{i}}>0$. In particular, for $\vec{q}=\vec{r}$, they cover the mixed-norm Besov-Lebesgue space ${\cal \dot{B}}^{-\frac{2}{p}}_{\vec{q},p}(\mathbb{R}^{d})$, for  $\vec{q}\in (2,\infty)^{d}$ and $2<p<\infty$ such that $\frac{-2}{p}=-1+\displaystyle{\sum_{i=1}^{d}}\frac{1}{q_{i}}$. 

In our previous work \cite{Leithold-Neves-Besov-mixed-Lebesgue}, we proved two results on global well-posedness for the Navier-Stokes system for small initial data in Besov type spaces based on 
mixed-Lebesgue spaces:

The first result is on the existence of a unique global mild solution, for small initial data in the Besov mixed-Lebesgue spaces 
\begin{equation*}
\dot{{\cal B}}^{\sigma}_{\vec{p},q}(\mathbb{R}^{d}),\,\,\mbox{where}\,\,\sigma=-1+\displaystyle{\sum_{i=1}^{d}}\frac{1}{p_{i}},\,\, \mbox{for}\,\, 1\leq q\leq \infty\,\, \mbox{and}\,\,\vec{p}\in [1,\infty]^{d}\,\,\mbox{such that}\,\,\displaystyle{\sum_{i=1}^{d}}\frac{1}{p_{i}}>0.
\end{equation*}
In particular, we need at least some $p_{i}\in [1,\infty)$. In fact, for $p_{1}=\dots=p_{d}=\infty$, 
the system \eqref{BQ} is ill-posed, see \cite{Bourgain-Pavlovic-2008}. Therefore, that work extended the results by Furioli, Lemarie and Terrano \cite{Furioli-Lemarie-Terrano-limite-2000}, Khai and Tri \cite{Khai-Tri-Sobolev-Lorentz-Navier-Stokes-2014} 
and by Bahouri, Chemin and Danchin (\cite{Bahouri-Chemin-Danchin-Fourier-PDE-2011}, Subsection 5.6.1). 
Here we also complement the previous and recent global well-posedness result in \cite{Phan-mixed-2020}.
Moreover, as a direct implication we can say that, for the next family of Besov type spaces $\dot{{\cal B}}^{-1+\frac{1}{p}}_{(\infty,\infty,p),\infty}(\mathbb{R}^{3})$, for $1\leq p\leq \infty$, 
our result provides a new optimal class for the well-posedness for $1\leq p< \infty$, since for $p=\infty$ the system is ill-posedness as proved in 
the well-known result \cite{Bourgain-Pavlovic-2008} by Bourgain and Pavlovi\'c. 

\medskip
The second result is on the existence of a unique global mild solution, for initial data in the Fourier-Besov mixed-Lebesgue spaces $\dot{{\cal FB}}^{s}_{\vec{p},q}(\mathbb{R}^{d})$,
\begin{equation*}
\mbox{where}\,\,s=-1+\displaystyle{\sum_{i=1}^{d}}\left(1-\frac{1}{p_{i}}\right),\,\, \mbox{for}\,\, 1\leq q\leq \infty\,\, \mbox{and}\,\,\vec{p}\in [1,\infty]^{d}\,\,\mbox{such that}\,\,\displaystyle{\sum_{i=1}^{d}}\left(1-\frac{1}{p_{i}}\right)>0,    
\end{equation*}
 which allow to partly extend previous results by Konieczny and Yoneda \cite{Konieczny-Yoneda-2011}, Iwabuchi and Takada \cite{Iwabuchi-Takada-2014}, Lei and Lin \cite{Lein-Lin-2011} and Fang, Han and Hieber \cite{Fang-Han-Hieber-Global-Navier-Stokes-Coriolis-Fourier-Besov-2015}.  
 Finally, as a particular case we have that, for the family of Fourier-Besov type spaces $\dot{{\cal FB}}^{-\frac{1}{p}}_{(1,1,p),\infty}(\mathbb{R}^{3})$, where $1\leq p\leq \infty$, 
 our result provide a new optimal class for the well-posedness for $1< p\leq \infty$, 
 since for $p=1$ the system \eqref{BQ}
 is ill-posed due to the well-known result by Iwabuchi and Takada \cite{Iwabuchi-Takada-2014}. 

%%%%%%%%%%%%%%%%%%%%
\subsection{Structural and analytical features of new function spaces}
%%%%%%%%%%%%%%%%%%%%

In the actual work we introduce two new functional spaces of distributions, namely, the homoge-
neous Besov mixed-Morrey spaces and Fourier-Besov mixed-Morrey spaces. We also establish some
basic properties on these spaces, which are described in the next two theorems. Notations are provided in Section \ref{preliminar}.

The next theorem describes the main properties related to Besov mixed-Morrey spaces.
\begin{teo}\label{properties-Besov-spaces-mixed-Morrey}

The following itens are valid:
\begin{itemize}
\item [(i.)] 
Let $\vec{q}\in [1,\infty)^{d}$, $\vec{\lambda}\in [0,1)^{d}$ and $l\in\mathbb{Z}$. There exists a constant $C$, independent on $l$, such that
\begin{equation}
 \parallel \dot{\Delta}_{l}f\parallel_{L^{\infty}(\mathbb{R}^{d})}\leq C2^{l\left(\sum_{i=1}^{d}\frac{1-\lambda_{i}}{q_{i}}\right)} \parallel \dot{\Delta}_{l}f\parallel_{{\cal M}_{\vec{q},\vec{\lambda}}(\mathbb{R}^{d})},
\end{equation}
for all $f\in {\cal S}'/{\cal P}$, such that, $\dot{\Delta}_{l}f\in {\cal M}_{\vec{q},\vec{\lambda}}(\mathbb{R}^{d})$.
\item[(ii.)] Let $\vec{q},\vec{r}\in [1,\infty)^{d}$, $\vec{\lambda},\vec{\mu}\in [0,1)^{d}$, $a\in [1,\infty]$ and $\sigma\in\mathbb{R}$.
\begin{itemize}
    \item [(ii.1)] For $\vec{1}\leq \vec{q}\leq \vec{r}<\vec{\infty}$ and $\frac{\lambda_{i}}{q_{i}}\geq \frac{\mu_{i}}{r_{i}}$, for $1\leq i\leq d$. If $\frac{1-\lambda_{i}}{q_{i}}=\frac{1-\mu_{i}}{r_{i}}$, for $1\leq i\leq d$, then it is valid the continuous inclusion
    \begin{equation}
     \dot{\cal N}^{\sigma}_{\vec{r},\vec{\mu},a}(\mathbb{R}^{d})\hookrightarrow \dot{\cal N}^{\sigma}_{\vec{q},\vec{\lambda},a}(\mathbb{R}^{d}).
    \end{equation}
    \item[(ii.2)] For $b\in[1,\infty]$, if $a\leq b$ then we have the continuous inclusion
    \begin{equation}
     \dot{\cal N}^{\sigma}_{\vec{q},\vec{\lambda},a}(\mathbb{R}^{d})\hookrightarrow \dot{\cal N}^{\sigma}_{\vec{q},\vec{\lambda},b}(\mathbb{R}^{d}).
     \end{equation}
\end{itemize}
\item [(iii.)] For $\vec{q}\in [1,\infty)^{d}$, $\vec{\lambda}\in [0,1)^{d}$, $r\in[1,\infty]$ and $s\in\mathbb{R}$ we have the continuous inclusion
\begin{equation}
     \dot{\cal N}^{s}_{\vec{q},\vec{\lambda},r}(\mathbb{R}^{d})\hookrightarrow \dot{\cal B}^{s-\left(\sum_{i=1}^{d}\frac{1-\lambda_{i}}{q_{i}}\right)}_{\infty,r}(\mathbb{R}^{d}).
\end{equation}
\item[(iv.)] For $\vec{q}\in [1,\infty)^{d}$, $\vec{\lambda}\in [0,1)^{d}$, $r\in[1,\infty]$, $s\in\mathbb{R}$ and $\theta\in(0,1)$ we have the continuous inclusion
    \begin{equation}
     \dot{\cal N}^{s}_{\vec{q},\vec{\lambda},r}(\mathbb{R}^{d})\hookrightarrow \dot{\cal N}^{s-\left(\sum_{i=1}^{d}\frac{1-\lambda_{i}}{q_{i}}\right)\left(1-\theta\right)}_{\vec{q}/\theta,\vec{\lambda},r}(\mathbb{R}^{d}).
     \end{equation}
\item[(v.)]  For $\vec{q}\in [1,\infty)^{d}$, $\vec{\lambda}\in [0,1)^{d}$, $r\in[1,\infty]$ and $s\in\mathbb{R}$, the space $\dot{\cal N}^{s}_{\vec{q},\vec{\lambda},r}(\mathbb{R}^{d})$ is a Banach space and we have the continuous inclusion
\begin{equation}
 \dot{\cal N}^{s}_{\vec{q},\vec{\lambda},r}(\mathbb{R}^{d})\hookrightarrow {\cal S}'/{\cal P}.   
\end{equation}
\item[(vi.)] Let us denote by $x=(x',x'')\in \mathbb{R}^{d_{2}}$, for $x'\in\mathbb{R}^{d_{1}}$ and $x''\in\mathbb{R}^{d_{2}-d_{1}}$, where $d_{1}<d_{2}$ are positive integers. Let us consider $s\in \mathbb{R}$, $\vec{1}\leq \vec{q}=(q_{1},\dots, q_{d_{1}},\dots, q_{d_{2}})<+\vec{\infty}$, $\vec{0}\leq \vec{\lambda}=(\lambda_{1},\dots,\lambda_{d_{1}})<\vec{1}$ and $\vec{0}\leq \vec{\beta}=(\beta_{1},\dots,\beta_{d_{2}})<\vec{1}$. Then, every element $u(x')$ in $\dot{{\cal N}}^{s}_{\vec{q},\vec{\lambda}}(\mathbb{R}^{d_{1}})$ is also an element in $\dot{{\cal N}}^{s}_{\vec{q},\vec{\beta}}(\mathbb{R}^{d_{2}})$, if the next condition
\begin{equation}
 \displaystyle{\sum_{i=1}^{d_{2}}}\frac{1-\beta_{i}}{q_{i}}=\displaystyle{\sum_{i=1}^{d_{1}}}\frac{1-\lambda_{i}}{q_{i}}   
\end{equation}
is valid.    
\end{itemize}   
\end{teo}
The next theorem describes the main properties related to Fourier-Besov mixed-Morrey spaces.

\begin{teo}\label{fourier-besov-properties}

The following itens are valid:
\begin{itemize}
\item [(i.)] 
Let $\vec{q}\in [1,\infty)^{d}$, $\vec{\lambda}\in [0,1)^{d}$ and $l\in\mathbb{Z}$. There exists a constant $C>0$, independent on $l$, such that
\begin{equation}
 \parallel \phi_{l} \hat{f}\parallel_{L^{1}(\mathbb{R}^{d})}\leq C 2^{l\left(\sum_{i=1}^{d}\left(1-\frac{1-\lambda_{i}}{q_{i}}\right)\right)} \parallel \phi_{l} \hat{f}\parallel_{{\cal M}_{\vec{q},\vec{\lambda}}(\mathbb{R}^{d})},
\end{equation}
for all $f\in {\cal S}'/{\cal P}$, such that, $\phi_{l} \hat{f}\in {\cal M}_{\vec{q},\vec{\lambda}}(\mathbb{R}^{d})$.
\item[(ii.) ] Let $\vec{1}\leq \vec{q}<\vec{r}<\vec{\infty}$ and $\vec{\lambda},\vec{\mu}\in[0,1)^{d}$ be parameters such that $\frac{\mu_{i}}{r_{i}}\leq \frac{\lambda_{i}}{q_{i}}$ and $\frac{1-\mu_{i}}{r_{i}}< \frac{1-\lambda_{i}}{q_{i}}$, for $1\leq i\leq d$, and $l\in\mathbb{Z}$. Then, there exists a positive constant $C$, independent on $l$, such that
\begin{equation}
 \parallel \phi_{l}\hat{f}\parallel_{{\cal M}_{\vec{q},\vec{\lambda}}(\mathbb{R}^{d})}\leq C2^{l\left(\sum_{i=1}^{d}\left(\frac{1-\lambda_{i}}{q_{i}}-\frac{1-\mu_{i}}{r_{i}}\right)\right)}\parallel \phi_{l}\hat{f}\parallel_{{\cal M}_{\vec{r},\vec{\mu}}(\mathbb{R}^{d})},   
\end{equation}
for each $f\in{\cal S}'/{\cal P}$, such that, $\phi_{l}\hat{f}\in{\cal M}_{\vec{r},\vec{\mu}}(\mathbb{R}^{d})$.
\item[(iii.)] Let $\vec{q},\vec{r}\in [1,\infty)^{d}$, $\vec{\lambda},\vec{\mu}\in [0,1)^{d}$, $a\in [1,\infty]$ and $\sigma\in\mathbb{R}$.
\begin{itemize}
    \item [(iii.1)] For $\vec{1}\leq \vec{q}\leq \vec{r}<\infty$ and $\frac{\lambda_{i}}{q_{i}}\geq \frac{\mu_{i}}{r_{i}}$, for $1\leq i\leq d$. If $\frac{1-\lambda_{i}}{q_{i}}=\frac{1-\mu_{i}}{r_{i}}$, for $1\leq i\leq d$,  then it is valid the continuous inclusion
    \begin{equation}
     \dot{\cal FN}^{\sigma}_{\vec{r},\vec{\mu},a}(\mathbb{R}^{d})\hookrightarrow \dot{\cal FN}^{\sigma}_{\vec{q},\vec{\lambda},a}(\mathbb{R}^{d}).
    \end{equation}
    \item[(iii.2)] For $b\in[1,\infty]$, if $a\leq b$ then we have the continuous inclusion
    \begin{equation}
     \dot{\cal FN}^{\sigma}_{\vec{q},\vec{\lambda},a}(\mathbb{R}^{d})\hookrightarrow \dot{\cal FN}^{\sigma}_{\vec{q},\vec{\lambda},b}(\mathbb{R}^{d}).
     \end{equation}
\end{itemize}
\item[(iv.)] Let $\vec{1}\leq \vec{q}<\vec{r}<\vec{\infty}$, $\vec{\lambda},\vec{\mu}\in[0,1)^{d}$, $a\in[1,\infty]$ and $s_{1},s_{2}\in\mathbb{R}$ be parameters such that $\frac{\mu_{i}}{r_{i}}\leq \frac{\lambda_{i}}{q_{i}}$ and $\frac{1-\mu_{i}}{r_{i}}< \frac{1-\lambda_{i}}{q_{i}}$, for $1\leq i\leq q$, and $s_{2}<s_{1}$. If $s_{2}+\sum_{i=1}^{d}\frac{1-\lambda_{i}}{q_{i}}=s_{1}+\sum_{i=1}^{d}\frac{1-\mu_{i}}{r_{i}}$ then
\begin{equation}
\dot{\cal FN}^{s_{1}}_{\vec{r},\vec{\mu},a}(\mathbb{R}^{d})\hookrightarrow \dot{\cal FN}^{s_{2}}_{\vec{q},\vec{\lambda},a}(\mathbb{R}^{d})
\end{equation}
\item [(v.)] For $\vec{q}\in [1,\infty)^{d}$, $\vec{\lambda}\in [0,1)^{d}$, $r\in[1,\infty]$ and $\sigma\in\mathbb{R}$ we have the next continuous inclusions
\begin{itemize}
\item [(v.1)]
\begin{equation}
\dot{\cal FN}^{\sigma}_{\vec{q},\vec{\lambda},r}(\mathbb{R}^{d})\hookrightarrow \dot{\cal FB}^{\sigma-\sum_{i=1}^{d}\left(1-\frac{1-\lambda_{i}}{q_{i}}\right)}_{1,r}(\mathbb{R}^{d}),\,\,\mbox{and}
\end{equation}
\item [(v.2)]
\begin{equation}
\dot{{\cal FB}}^{\sigma}_{1,r}(\mathbb{R}^{d})\hookrightarrow \dot{{\cal B}}^{\sigma}_{\infty,r}(\mathbb{R}^{d}).    
\end{equation}
\end{itemize}

\item[(vi.)] For $\vec{q}\in [1,\infty)^{d}$, $\vec{\lambda}\in [0,1)^{d}$, $r\in[1,\infty]$, $s\in\mathbb{R}$ and $\theta\in(0,1)$ we have the continuous inclusion
    \begin{equation}
     \dot{\cal FN}^{s}_{\vec{q},\vec{\lambda},r}(\mathbb{R}^{d})\hookrightarrow \dot{\cal FN}^{s-\sum_{i=1}^{d}(\frac{1-\lambda_{i}}{q_{i}})(1-\theta)}_{\vec{q}/\theta,\vec{\lambda},r}(\mathbb{R}^{d}).
     \end{equation}
\item[(vii.)]  For $\vec{q}\in [1,\infty)^{d}$, $\vec{\lambda}\in [0,1)^{d}$, $r\in[1,\infty]$ and $s\in\mathbb{R}$, the space $\dot{\cal FN}^{s}_{\vec{q},\vec{\lambda},r}(\mathbb{R}^{d})$ is a Banach space and we have the continuous inclusion
\begin{equation}
 \dot{\cal FN}^{s}_{\vec{q},\vec{\lambda},r}(\mathbb{R}^{d})\hookrightarrow {\cal S}'/{\cal P}.   
\end{equation}
\end{itemize}   
\end{teo}

%%%%%%%%%%%%%%%%%%%
\subsection{On an application} 
%%%%%%%%%%%%%%%%%%%

The above results are used to establish the next two theorems of global well-posedness, for small initial data, for the $d$-dimensional incompressible Navier-Stokes system, 
%In this work we deal with the $d$-dimensional incompressible Navier-Stokes equations, that is to say, 
\begin{equation}
\label{BQ}
%(NS)^{(\nu)}\,\,\,
\left\{
\begin{aligned}
&  \partial_{t} u(x,t) + (u(x,t) \cdot \nabla) u(x,t) 
+ \nabla P(x,t) = \nu \Delta u(x,t),
\\[5pt]
&  \mbox{div}\,u(x,t)= 0,\quad \mbox{for}\,\,x=(x_1, \dots, x_d) \in \mathbb{R}^{d}, \quad t> 0,
\\[5pt]
&  u(x,0)=u_{0}(x),\,\,\mbox{for}\,x\in\mathbb{R}^{d},
\end{aligned}
\right.  
\end{equation}
%%%%%%%%%%%%%%%%%%%%%%%%%%%%%%%%%%%%%%%%%
where $u=u(x,t)=(u^{1}(x,t), \dots, u^{d}(x,t))$ represents the velocity vector field, the scalar function $P=P(x,t)$ denotes the pressure, 
the parameter $\nu> 0$ represents the kinematic viscosity, and $u_{0}$ is a given divergence free initial velocity field. 

A direct application of Besov mixed-Morrey spaces by mean of a Bony's paraproduct approach is given in the next

\begin{teo}\label{teorema-besov-morrey-misto}
Let $\vec{q}\in[1,\infty)^{d}$, $\vec{\lambda}\in[0,1)^{d}$ and $r\in[1,\infty]$ be parameters such that $\displaystyle{\sum_{i=1}^{d}}\frac{1-\lambda_{i}}{q_i}>0$, $\sigma=-1+\displaystyle{\sum_{i=1}^{d}}\frac{1-\lambda_{i}}{q_i}$ and
$Z={\cal L}^{\infty}\left(I;{\dot{{\cal N}}^{\sigma}_{\vec{q},\vec{\lambda},r}(\mathbb{R}^{d})}\right)\bigcap {\cal L}^{1}\left(I;{\dot{{\cal N}}^{\sigma+2}_{\vec{q},\vec{\lambda},r}(\mathbb{R}^{d})}\right)$, where $I=(0,\infty)$. Then,
\begin{itemize}
\item [(i.)] There exist $\delta>0$ and $D_0$ such that, for small enough initial data $\parallel u_{0}\parallel_{\dot{{\cal N}}^{\sigma}_{\vec{q},\vec{\lambda},r}(\mathbb{R}^{d})}\leq \delta$ there is a unique global mild solution $u$ in the ball $\{u\in Z;\parallel u\parallel_{Z}\leq \delta D_0\}$. Moreover, the map solution depends continuously on the initial data. 
\item [(ii.)] The obtained solution is weakly continuous from $[0,\infty)$ to $S'(\mathbb{R}^{d})$. 
\end{itemize}
\end{teo}

\begin{obs}
 First, let us observe that the condition $\sum_{i=1}^{d}\frac{1-\lambda_{i}}{q_{i}}>0$ is true for the considered parameters $\vec{q}$ and $\vec{\lambda}$. Moreover, this theorem partially extends the corresponding result in Besov-Morrey spaces (\cite{Kozono-Yamazaki-semilinear-1994}, Theorem 3).
\end{obs}

\begin{obs}\label{initial-data-example}
 In this remark we consider two examples of initial data for which the above theorem applies. 
 \begin{itemize}
\item [i.)] Let us consider $\vec{1}\leq \vec{q}<+\vec{\infty}$ and $\vec{0}<\vec{\lambda}<\vec{1}$. Take 
\begin{equation*}
f(x)=\prod_{j=1}^{d}\mid x_{j}\mid^{-\left(\frac{1-\lambda_{j}}{q_{j}}\right)}\in{\cal M}_{\vec{q},\vec{\lambda}}(\mathbb{R}^{d}).
\end{equation*}

Since (see item (i.) from Proposition \ref{auxiliar-propo-inclusion})

\begin{equation*}
\dot{{\cal N}}^{0}_{\vec{q},\vec{\lambda},1}(\mathbb{R}^{d})\subset {\cal M}_{\vec{q},\vec{\lambda}}(\mathbb{R}^{d})\subset \dot{{\cal N}}^{0}_{\vec{q},\vec{\lambda},\infty}(\mathbb{R}^{d}),
\end{equation*}
then $g=(-\Delta)^{\frac{\sigma}{2}}f\in \dot{{\cal N}}^{\sigma}_{\vec{q},\vec{\lambda},\infty}(\mathbb{R}^{d})$, for $\sigma=-1+\displaystyle{\sum_{i=1}^{d}\frac{1-\lambda_{i}}{q_{i}}}$.
\item [ii.)] Let us consider $\vec{1}\leq \vec{q}=(q_{1},\dots,q_{d})<+\vec{\infty}$, $0<\lambda_{1}<1$ and $\vec{0}\leq \vec{\beta}<\vec{1}$ satisfying the condition (\ref{condition-Morrey}). Then, since 
\begin{equation*}
\mid x_{1}\mid^{-\left(\frac{1-\lambda_{1}}{q_{1}}\right)}\in{\cal M}_{q_{1},\lambda_{1}}(\mathbb{R}),
\end{equation*}
we have by the condition (\ref{condition-Morrey}) that
\begin{equation*}
\mid x_{1}\mid^{-\left(\frac{1-\lambda_{1}}{q_{1}}\right)}\in{\cal M}_{\vec{q},\vec{\beta}}(\mathbb{R}^{d})\subset\dot{{\cal N}}^{0}_{\vec{q},\vec{\beta},\infty}(\mathbb{R}^{d}).    
\end{equation*}
Let us take $g=(-\Delta)^{\frac{s}{2}}\left(\mid x_{1}\mid^{-\left(\frac{1-\lambda_{1}}{q_{1}}\right)}\right)\in \dot{{\cal N}}^{s}_{\vec{q},\vec{\beta},\infty}(\mathbb{R}^{d})$, for $s=-1+\sum_{i=1}^{d}\frac{1-\beta_{i}}{q_{i}}$. 
 \end{itemize}
 For both the itens, the initial velocity $u_{0}=\varepsilon(0,\dots,0,g)$, for small enough $\varepsilon>0$, could be considered in order to apply the above Theorem \ref{teorema-besov-morrey-misto}.
\end{obs}

Another direct application of Fourier-Besov mixed-Morrey spaces by mean of a Bony's paraproduct approach is given in the next

\begin{teo}
Let $\vec{q}\in[1,\infty)^{d}$, $\vec{\lambda}\in[0,1)^{d}$ and $r\in[1,\infty]$ be parameters such that $\displaystyle{\sum_{i=1}^{d}}\left(1-\frac{1-\lambda_{i}}{q_i}\right)>0$, $s=-1+\displaystyle{\sum_{i=1}^{d}}\left(1-\frac{1-\lambda_{i}}{q_i}\right)$ and $Z={\cal L}^{\infty}\left(I;{\dot{{\cal FN}}^{s}_{\vec{q},\vec{\lambda},r}(\mathbb{R}^{d})}\right)\bigcap {\cal L}^{1}\left(I;{\dot{{\cal FN}}^{s+2}_{\vec{q},\vec{\lambda},r}(\mathbb{R}^{d})}\right)$, where $I=(0,\infty)$. Then,
\begin{itemize}
\item [(i.)] There exist $\delta>0$ and $D_0$ such that, for small enough initial data $\parallel u_{0}\parallel_{\dot{{\cal FN}}^{s}_{\vec{q},\vec{\lambda},r}(\mathbb{R}^{d})}\leq \delta$ there is a unique global mild solution $u$ in the ball $\{u\in Z;\parallel u\parallel_{Z}\leq \delta D_0\}$. Moreover, the map solution depends continuously on the initial data. 
\item [(ii.)] The obtained solution is weakly continuous from $[0,\infty)$ to $S'(\mathbb{R}^{d})$. 
\end{itemize}
\end{teo}
\begin{obs}
Particular cases from the condition $\displaystyle{\sum_{i=1}^{d}}\left(1-\frac{1-\lambda_{i}}{q_i}\right)>0$ could be considered, for instance, if $\vec{q}=\vec{1}$ then it is enough to take $\displaystyle{\sum_{i=1}^{d}}\lambda_{i}>0$, and if $q_{i}=1$, for $i=1,\dots,d-1$, then it is enough to consider 
\begin{equation*}
 q_{d}>\frac{1-\lambda_{d}}{1+\displaystyle{\sum_{i=1}^{d-1}}\lambda_{i}}.
\end{equation*}
\end{obs}

\begin{obs}
In order to complement the context of these main theorems we also mention the works on anisotropic Besov type spaces, by considering initial data in ${\cal B}^{0,\frac{1}{2}}(\mathbb{R}^{3})$ \cite{Paicu-anisotrope-2005} by Paicu, or in ${\cal B}^{-\frac{1}{2},\frac{1}{2}}_{4}(\mathbb{R}^{3})$ \cite{Paicu-Zhang-global-2011} by Paicu and Zhang.
\end{obs}

%%%%%%%%%%%%%%%%%%%%%%%%%%%%%%%%%%%%%%%%%%%%%%
The employed spaces and techniques in this work could be adapted to study some another incompressible models on fluid dynamics, as the surface quasi-geostrophic equation, the 3D-micropolar fluid system, the 3D-Keller-Segel-Navier-Stokes equations, the 3D-attraction-repulsion chemotaxis-fluid system, the 3D-Hall-magnetohydrodynamic system or the 3D-nematic liquid crystal equations. In fact there exist several results for these equations in the context of Besov spaces or Besov-Morrey spaces.
%%%%%%%%%%%%%%%%%%%
\subsection{Related isotropic functional spaces}
In order to complement the context of the above results, we mention some results on existence, uniqueness or stability for the Navier-Stokes system, in different kind of spaces.

For instance, in Morrey spaces \cite{Kato-strong-Morrey-1992}, ${\cal M}_{q,\lambda}(\mathbb{R}^{d})$, with $\lambda=d-q$ and $1\leq q<\infty$, and in the weak Morrey space \cite{Miao-Yuan-weak-Morrey-2007}, ${\cal M}^{\ast}_{p,\lambda}(\mathbb{R}^{d})$, for $1<p\leq d$ and $\lambda=d-p$, in the space 
\cite{Planchon-NS-1996}, $H^{s}(\mathbb{R}^{3})$, for $s>0$, satisfying a smallness condition in the Besov space $\dot{{\cal B}}^{-\frac{1}{4}}_{4,\infty}(\mathbb{R}^{3})$, and in the space \cite{Koch-Tataru-2001}, $BMO^{-1}(\mathbb{R}^{3})$, which is up to the present, one of the maximal classes for this kind of results.

In the modern textbook (\cite{Bahouri-Chemin-Danchin-Fourier-PDE-2011}, 
Subsection 5.6.1) Bahouri, Chemin and Danchin revisited several results in Besov spaces and established the global well-posedness for initial data 
$u_{0}\in\dot{{\cal B}}^{-1+\frac{3}{p}}_{p,\infty}(\mathbb{R}^{3})$, for $1\leq p<\infty$. 
A recently preprint \cite{Cheskidov-Eguchi-NS-Koch-2025} explore the existence of a unique global solution for initial data in $L^{2}(\mathbb{R}^{3})$, by considering a smallness condition for the high-frequency in the space $BMO^{-1}(\mathbb{R}^{3})$ and for the low-frequency in the space $\dot{\cal B}^{-1}_{\infty,\infty}(\mathbb{R}^{3})$. For an ill-posedness resut in a Besov space, Bourgain and Pavlovi\'c \cite{Bourgain-Pavlovic-2008} proved that the solution map for the Navier-Stokes equations is 
not continuous in $\dot{{\cal B}}^{-1}_{\infty,\infty}(\mathbb{R}^{3})$ at the origin, which implies the ill-posedness, for initial 
data $u_{0}\in\dot{{\cal B}}^{-1}_{\infty,\infty}(\mathbb{R}^{3})$. 

Moreover, in the Fourier-Besov space \cite{Iwabuchi-Takada-2014}, $\dot{{\cal FB}}^{-1}_{1,2}(\mathbb{R}^{3})$, whose smallness condition does not depend on the Coriolis parameter, and in the space \cite{Fang-Han-Hieber-Global-Navier-Stokes-Coriolis-Fourier-Besov-2015}, $\dot{{\cal FB}}^{2-\frac{3}{p}}_{p,r}(\mathbb{R}^{3})$, for $1<p\leq \infty$ and $1\leq r\leq\infty$. For an ill-posedness result in Fourier-Besov spaces  $\dot{{\cal FB}}^{-1}_{1,q}(\mathbb{R}^{3})$, for any $2< q\leq \infty$, see the work \cite{Iwabuchi-Takada-2014}, in Besov-Morrey spaces \cite{Kozono-Yamazaki-semilinear-1994}, ${\cal N}^{-1+\frac{d-\lambda}{q}}_{q,\lambda,\infty}(\mathbb{R}^{d})$, for $1\leq q<\infty$ and $\lambda\in[0,d)$ such that $d-q<\lambda$ and in Fourier-Besov-Morrey spaces \cite{Almeida-Ferreira-Lima-uniform-2017}, $\dot{{\cal FN}}^{2-\frac{3-\mu}{q}}_{q,\mu,\infty}(\mathbb{R}^{3})$, for $0\leq \mu<3$ and $1\leq q<\infty$ ($\mu\neq 0$ if $q=1$). These last two spaces are each the counterpart of the another, that is, as the first one analyse the system considering a localization in Morrey spaces, the second one set the localization in Fourier-variables based on Morrey spaces.
%%%%%%%%%%%%%%%%%%%%%%%%%

In a recent work, Ferreira and P\'erez-L\'opez \cite{Ferreira-Perez-Lopez-BWherz-2018} proved the existence of a unique global solution of the $(NS)$-system, for small initial data in Besov-weak-Herz spaces $\dot{{\cal B}}W\dot{K}^{\alpha,\alpha+\frac{d}{p}-1}_{p,q,\infty}(\mathbb{R}^{d})$, where $1 \leq q\leq \infty$, $\frac{d}{2}<p<\infty$, and $0\leq \alpha<\min\{1-\frac{d}{2p},\frac{d}{2p}\}$. As indicated in that work, these spaces satisfy the next continuous inclusions
\begin{equation}
\begin{split}
 L^{d}(\mathbb{R}^{d})\hookrightarrow L^{d,\infty}(\mathbb{R}^{d})\hookrightarrow W\dot{K}^{0}_{d,\infty}(\mathbb{R}^{d})\hookrightarrow \dot{{\cal B}}W\dot{K}^{0,0}_{d,\infty,\infty}(\mathbb{R}^{d})\,\,\mbox{and}\\
 \dot{H}^{-1+\frac{d}{2}}(\mathbb{R}^{d})\hookrightarrow L^{d,\infty}(\mathbb{R}^{d})\hookrightarrow \dot{{\cal B}}^{-1+\frac{d}{p}}_{p,\infty}(\mathbb{R}^{d})\hookrightarrow \dot{{\cal B}}W\dot{K}^{0,-1+\frac{d}{p}}_{p,\infty,\infty}(\mathbb{R}^{d})\,\,\mbox{for}\,\,p\geq d. 
\end{split}    
\end{equation}
Moreover, this new class is to the best knowledge of the authors one of the maximal class of well-posedness for the system.

Finally, we can mention the work \cite{Abidin-Chen-variable-FBM-2021}, for the generalized Navier-Stokes system, in the context of variable Fourier-Besov-Morrey spaces  $\dot{{\cal FN}}^{4-2\alpha-\frac{3}{p(\cdot)}}_{p(\cdot),h(\cdot),q}(\mathbb{R}^{3})$, for $\frac{1}{2}<\alpha\leq 1$, $p(\cdot),h(\cdot)\in C^{log}(\mathbb{R}^{3})\cap{\cal P}_{0}(\mathbb{R}^{3})$ such that $2\leq p(\cdot)\leq\frac{6}{5-4\alpha}$ and $p(\cdot)\leq h(\cdot)<\infty$, and $1\leq q<\frac{3}{2\alpha-1}$. 
%%%%%%%%%%%%%%%%%%%%%%%%%%%%%%%%%%%%%%%%%%%%%%%%%%%%%%%%%%%%%%%%%%%%%%%%%%

%%%%%%%%%%%%%%%%%%%%%%%%%%%%%%%%%%%%%%%%%%%%%%%%%%%%%%%%%%%%%%

%\medskip
%From all of the previous cited works, one may say that, the use of Besov spaces (also Fourier-Besov spaces) has led to various advancements 
%in existence, uniqueness, and regularity results for weak 
%and mild solutions, particularly in spaces that allow low regularity initial data. Moreover, Besov spaces often allow us to establish 
%global existence of solutions in the critical space (where the scaling invariance holds) for small initial data, which is 
%a crucial aspect in fluid mechanics and turbulence theory.

%\medskip
%Albeit, the theory of  Besov type spaces, for instance, anisotropic Besov spaces, Besov spaces based on mixed-Lebesgue spaces, %
%Besov spaces based on mixed-Lorentz spaces, etc. are under development. These spaces are extensions of classical Besov spaces,
%which are designed to handle directional or spatially dependent irregularities. Indeed, in many physical and mathematical problems, 
%the solution may exhibit different smoothness properties in different spatial directions. Then, mixed-norm Besov type spaces allow 
%separate control over regularity along each spatial direction. This flexibility is crucial for fluid dynamics problems, 
%especially when flows have directional structures, like in Navier-Stokes equations. 
%
%In particular, we cite below two recently works in this direction.

%%%%%%%%%%%%%%%%%%%%%%%%%%%

\bigskip
The structure of this paper is the following:

\medskip
In Section \ref{preliminar}, we recall Morrey-spaces, mixed-Lebesgue spaces, mixed-Morrey spaces, Besov mixed-Lebesgue spaces, Fourier-Besov mixed-Lebesgue spaces. In Section \ref{section-BMMorrey} we define Besov mixed-Morrey spaces and prove some basic properties on these spaces.
In Section \ref{section-FBMMorrey} we define Fourier-Besov mixed-Morrey spaces and prove some of its basic properties.
In the last Section \ref{Section-NS-system} we provide a direct application, in the context of the Navier-Stokes equations, of the introduced Besov mixed-Morrey spaces and Fourier-Besov mixed-Morrey spaces. First, we establish the linear and bilinear estimates for the mild formulation and finally we prove the last two main theorems.

%%%%%%%%%%%%%%%%%%%%%%%%%%%%%%%%%
%%%%%%%%%%%%%%%%%%%%%%%%%%%%%%%%
%%%%%%%%%%%%%%%%%%%%%%%%%%%%%%%%
%%%%%%%%%%%%%%%%%%%%%%%%%%%%%%%%%%%
%%%%%%%%%%%%%%%%%%%%%%%%%%%%%%%%%%%%%%
%%%%%%%%%%%%%%%%%%%%%%%%%%%%%%
\section{Preliminaries}\label{preliminar}

In this section, we fix the notations and collect some preliminary results. 
First, let $\Omega \subset \mathbb{R}^d$ be an open set. 
We denote by $dx, d\xi$, etc. 
the Lebesgue measure on $\Omega$ and by $L^p(\Omega)$, $p \in [1,+\infty)$, the set of $p$-summable functions 
with respect to the Lebesgue measure, also $L^\infty(\Omega)$ is 
the set of measurable functions in which its absolute value has the essential supremum finite. 
As is customary, the symbols $\mathcal{S}(\mathbb{R}^d)$ denotes the Schwartz class, and
$\mathcal{S}^{\prime}(\mathbb{R}^d)$ denotes the set of tempered distributions. We denote by $\mathcal{F} \varphi(\xi) \equiv \widehat{\varphi}(\xi)\equiv\varphi^{\wedge}(\xi)$
the Fourier Transform of $\varphi$, which is an
isometry in $L^2(\mathbb{R}^d)$. Moreover, we denote by ${\cal D}(G)$ the set of infinitely differentiable real functions whose support is contained in the subset $G$ of $\mathbb{R}^{d}$.

\medskip
Let us recall the Dyadic Partition of Unity from \cite{Bahouri-Chemin-Danchin-Fourier-PDE-2011}.  
\begin{lema}
Let ${\cal A}=\{\eta\in \mathbb{R}^{d};\frac{3}{4}\leq \mid \eta\mid\leq \frac{8}{3}\}$ be a given annulus. There exist radial functions $\psi$ and $\phi$, valued in the interval $[0,1]$, belonging respectively to ${\cal D}(B(0,\frac{4}{3}))$ and ${\cal D}({\cal A})$, and such that
\begin{equation}
\begin{split}
 \forall \eta\in\mathbb{R}^{d},\,&\psi(\eta)+\displaystyle{\sum_{l\geq 0}\phi(2^{-l}\eta)}=1, \\
 \forall \eta\in\mathbb{R}^{d}-\{0\},\,\,&\displaystyle{\sum_{l\in\mathbb{Z}}\phi(2^{-l}\eta)}=1
\end{split}   
\end{equation}
and satisfy the following intersection properties for their supports
\begin{equation}
 \begin{split}
 \mid l-l'\mid\geq 2&\implies \mbox{supp}\,\phi(2^{-l}\cdot)\cap \mbox{supp}\,\phi(2^{-l'}\cdot)=\emptyset,\\
 l\geq 1&\implies \mbox{supp}\,\psi\, \cap\mbox{supp}\,\phi(2^{-l}\cdot)=\emptyset.
 \end{split}   
\end{equation}
\end{lema}

\begin{obs}
The function $\tilde{\psi}(\eta)=\psi(\eta)+\phi(\eta)$ is equal to one in some neighborhood of $B(0,\frac{4}{3})$. Moreover the function $\tilde{\phi}_{l}=\phi_{l-1}+\phi_{l}+\phi_{l+1}$, where $\phi_{l}(\eta)=\phi(2^{-l}\eta)$, is equal to one in some neighborhood of the annulus ${\cal A}_{l}=2^{l}{\cal A}=\{\eta\in \mathbb{R}^{d};\frac{3}{4}\cdot2^{l}\leq \mid \eta\mid\leq 2^{l}\cdot\frac{8}{3}\}$.
\end{obs}

\begin{obs}
Since $\sum_{l<0}\phi(2^{-l}\eta)=\psi(\eta)$, for all $\eta\in\mathbb{R}^{d}-\{0\}$, and by definition of the operator $\dot{S}_{l}g=\psi(2^{-l}D)g={\cal F}^{-1}[\psi(2^{-l}\eta)\hat{g}]$, it is valid that
\begin{equation}
  \dot{S}_{l}g=\displaystyle{\sum_{l'\leq l-1}}\dot{\Delta}_{l'}g, \,\,\mbox{for each}\, l\in\mathbb{Z},
\end{equation}
where the operator $\dot{\Delta}_{l}$ is given by $\dot{\Delta}_{l}g={\cal F}^{-1}[\phi_{l}\hat{g}]$. 
\end{obs}

The main properties of these operators are the next, so called, almost orthogonal ones, that is to say
\begin{equation*}
\begin{split}
\dot{\Delta}_l \dot{\Delta}_{l'}g&=0,\,\,\mbox{if}\,\,\mid l-l'\mid\geq 2,\,\,\mbox{and}\\
\dot{\Delta}_l\left(\dot{S}_{l'-1}f\dot{\Delta}_{l'}g \right)&=0, \,\,\mbox{if}\,\,\mid l-l'\mid\geq 5.
\end{split}
\end{equation*}

\subsection{Besov mixed-Lebesgue spaces}
In this section we recall the definition of Besov mixed-Lebesgue spaces and its main properties as the Bernstein-type lemmas, embedding inclusions and completeness. 

We first recall the mixed-Lebesgue spaces and its main properties.
 \begin{defi}(Mixed-Lebesgue spaces)
 Let us consider $\vec{p}=(p_{1},\dots, p_{d})\in [1,\infty]^{d}$. The Mixed-Lebesgue spaces $L^{\vec{p}}(\mathbb{R}^{d})$ is the set of measurable functions $f$ defined on $\mathbb{R}^{d}$ such that 
 \begin{equation*}
   \parallel f\parallel_{L^{\vec{p}}(\mathbb{R}^{d})}=\left(
   \displaystyle{\int_{\mathbb{R}}}\left( \displaystyle{\int_{\mathbb{R}}}\dots \left(\displaystyle{\int_{\mathbb{R}}}\mid f(x_{1},\dots, x_{d})\mid^{p_{1}}dx_{1}\right)^{\frac{1}{p_{1}}}\dots dx_{d-1}\right)^{\frac{p_{d}}{p_{d-1}}} dx_{d}\right)^{\frac{1}{p_{d}}}\,\,\mbox{is finite}.
 \end{equation*}
 Here we consider integrals if $\vec{p}\in [1,\infty)^{d}$ and the $i-$th essential sup norm if some $p_{i}=\infty$. 
 \end{defi}
Next lemma summarize some main properties on mixed-Lebesgue spaces.
\begin{lema} 
 \item [(i.)] For $\vec{p}_{1},\vec{p}_{2},\vec{p}_{3}\in [1,\infty]^{d}$, such that $\frac{1}{\vec{p}_{3}}=\frac{1}{\vec{p}_{1}}+\frac{1}{\vec{p}_{2}}$ we have
\begin{equation*}
 \parallel g_1\cdot g_2\parallel_{L^{\vec{p}_{3}}(\mathbb{R}^{d})}\leq    \parallel g_1\parallel_{L^{\vec{p}_{1}}(\mathbb{R}^{d})}\cdot
 \parallel g_2\parallel_{L^{\vec{p}_{2}}(\mathbb{R}^{d})},
\end{equation*}
for all $g_{1}\in L^{\vec{p_{1}}}(\mathbb{R}^{d})$ and $g_{2}\in L^{\vec{p_{2}}}(\mathbb{R}^{d})$.
\item  [(ii.)] For $\vec{p}\in[1,\infty]^{d}$ we have 
\begin{equation*}
 \parallel \phi\ast g\parallel_{L^{\vec{p}}(\mathbb{R}^{d})}\leq \parallel \phi\parallel_{L^{1}(\mathbb{R}^{d})}\parallel g\parallel_{L^{\vec{p}}(\mathbb{R}^{d})},
 \end{equation*}
 for all $\phi \in L^{1}(\mathbb{R}^{d})$ and $g\in L^{\vec{p}}(\mathbb{R}^{d})$.
 \item [(iii.)] The normed space $\left(L^{\vec{p}}(\mathbb{R}^{d}),\parallel \cdot\parallel_{L^{\vec{p}}(\mathbb{R}^{d})}\right)$ is complete.
\end{lema}
Here, we denote $\frac{1}{\vec{p}}=(\frac{1}{p_{1}},\dots,\frac{1}{p_{d}})$, where $\vec{p}=(p_{1},\dots,p_{d})\in[1,\infty]^{d}$.
For further properties on these spaces see the work \cite{Benek-Panzone-mixed-Lp-1961}, by A. Benedek and R. Panzone. 
 
\bigskip
In the modern definition, mixed-norm Besov spaces were introduced recently in the work \cite{Cleanthous-Georgiadis-Nielsen-mixed-norm-2016}, 
where they established some main properties as well as the discrete decomposition of homogeneous mixed-norm Besov spaces. 
To see another related studies we address the reader to \cite{Cleanthous-Georgiadis-Nielsen-discrete-2017} and \cite{Georgiadis-Nielsen-Pseudo-2016}, also 
for usual Besov spaces see \cite{Sawano-Besov-2018} and \cite{Bahouri-Chemin-Danchin-Fourier-PDE-2011}.
Before we recall the Besov mixed-Lebesgue spaces, we denote by ${\cal S}'/{\cal P}$ the space of tempered distributions modulo polynomials on $\mathbb{R}^{d}$. 

The next spaces were considered to study the Navier-Stokes system in \cite{Leithold-Neves-Besov-mixed-Lebesgue}.
\begin{defi}(Besov mixed-Lebesgue spaces)
\label{DefBesovType} 
Let us consider the three parameters 
$$
    \text{$\vec{p}=(p_1,\dots,p_d)\in [1,\infty]^{d}$, $q\in [1,\infty]$ and $\sigma\in\mathbb{R}$.}
$$    
The $d$-dimensional homogeneous Besov mixed-Lebesgue spaces are defined as
\begin{itemize}
\item [(i.)] For $\vec{p}\in [1,\infty]^{d}$, $q\in[1,\infty)$ and $\sigma\in\mathbb{R}$,
\begin{equation}
\dot{{\cal B}}^{\sigma}_{\vec{p},q}(\mathbb{R}^{d})=\lbrace
g\in {\cal S}'/{\cal P};\parallel g\parallel_{\dot{{\cal B}}^{\sigma}_{\vec{p},q}(\mathbb{R}^{d})}=
\left(\displaystyle{\sum_{l\in\mathbb{Z}}2^{l\sigma q}\parallel \dot{\Delta}_{l}g\parallel^{q}_{L^{\vec{p}}(\mathbb{R}^{d})}}
\right)^{\frac{1}{q}}\,\,\mbox{is finite}
\rbrace.
\end{equation}
\item [(ii.)]  For $\vec{p}\in [1,\infty]^{d}$, $q=\infty$ and $\sigma\in\mathbb{R}$,
\begin{equation}
\dot{{\cal B}}^{\sigma}_{\vec{p},\infty}(\mathbb{R}^{d})=\lbrace
g\in {\cal S}'/{\cal P};\parallel g\parallel_{\dot{{\cal B}}^{\sigma}_{\vec{p},\infty}(\mathbb{R}^{d})}=\displaystyle{\sup_{l\in\mathbb{Z}} 2^{l\sigma}\parallel \dot{\Delta}_{l}g\parallel_{L^{\vec{p}}(\mathbb{R}^{d})}}\,\,\mbox{is finite}
\rbrace.
\end{equation}
\end{itemize}
We also called these spaces as the $d$-dimensional homogeneous mixed-norm Besov-Lebesgue spaces or the $d$-dimensional homogeneous Besov spaces based on mixed-Lebesgue spaces.  
\end{defi}
\begin{lema}{(Main properties related to Besov mixed-Lebesgue spaces)}\label{properties-Besov-spaces-mixed}
\begin{itemize}
\item [(i.)]  Let ${\cal A}=\{\eta\in \mathbb{R}^{d}; r_{1}\leq \mid\eta\mid\leq r_{2}\}$, for some $0<r_{1}<r_{2}$, be a given annulus and 
$$
\text{$B=\{\eta\in \mathbb{R}^{d}; \mid \eta\mid\leq r\}$, for some $r>0$, be a closed ball.} 
$$
Then, there exist a positive constant $C$, such that, for any nonnegative integer $k$ and a couple 
$\vec{p},\vec{q}\in\left([1,\infty]^{d}\right)^{2}$, with $\vec{1}\leq\vec{p}\leq \vec{q}$, it follows for any function $u\in L^{\vec{p}}(\mathbb{R}^{d})$ that:
\begin{equation*}
 \begin{split}
   &(i.1.)\,\,\supp \hat{u}\subset \lambda B\implies\parallel D^{k}u\parallel_{L^{\vec{q}}(\mathbb{R}^{d})}=\displaystyle{\sup_{\mid \alpha\mid=k}}\parallel \partial^{\alpha}u\parallel_{L^{\vec{q}}(\mathbb{R}^{d})}\leq C^{k}\lambda^{k_{0}}\parallel u\parallel_{L^{\vec{p}}(\mathbb{R}^{d})},\,\,\mbox{and}\\
   &(i.2.)\,\, \supp\hat{u}\subset \lambda {\cal A}\implies C^{-k-1}\lambda^{k}\parallel u\parallel_{L^{\vec{p}}(\mathbb{R}^{d})}\leq \parallel D^{k}u\parallel_{L^{\vec
   p}(\mathbb{R}^{d})}\leq C^{k+1}\lambda^{k}\parallel u\parallel_{L^{\vec{p}}(\mathbb{R}^{d})},  
 \end{split}   
\end{equation*}
where $k_{0}=k+\displaystyle{\sum_{i=1}^{d}}\left(
   \frac{1}{p_i}-\frac{1}{q_i}\right)$.
\item[(ii.)] Let us consider $\vec{p},\vec{q}\in[1,\infty]^{d}$ such that $\vec{1}\leq\vec{p}\leq\vec{q}$ and $1\leq a_{1}\leq a_{2}\leq \infty$. For $\sigma_{2}\leq\sigma_{1}$ 
with $\sigma_{1}+\displaystyle{\sum_{i=1}^{d}}\frac{1}{q_{i}}=\sigma_{2}+\displaystyle{\sum_{i=1}^{d}}\frac{1}{p_{i}}$, we have
\begin{equation*}
 \dot{{\cal B}}^{\sigma_{1}}_{\vec{p},a_{1}}(\mathbb{R}^{d})\hookrightarrow
 \dot{{\cal B}}^{\sigma_{2}}_{\vec{q},a_{2}}(\mathbb{R}^{d}),
\end{equation*}
that is, the above inclusion is continuous. Here, $\vec{1}=(1,\dots,1)$ and $\vec{p}\leq\vec{q}$ means $p_{1}\leq q_{1},\dots, p_{d}\leq q_{d}$ if $\vec{p}=(p_{1},\dots,p_{d})$ and $\vec{q}=(q_{1},\dots,q_{d})$.
%{\bf nao esta definido o que seja $\vec{p}\leq \vec{q}$} 
\item [(iii.)] For $\vec{p}\in[1,\infty]^{d}$, $q\in[1,\infty]$ and $\sigma\in\mathbb{R}$, the Besov mixed-Lebesgue space $\dot{{\cal B}}^{\sigma}_{\vec{p},q}(\mathbb{R}^{d})$ is a Banach space endowed with the norm $\parallel \cdot\parallel_{\dot{{\cal B}}^{\sigma}_{\vec{p},q}(\mathbb{R}^{d})}$.
\end{itemize}   
\end{lema}
%%%%%%%%%%%%%%%%%%%%%%%%%%%%%%%%%%%%%%%%%%%%%%%%%%%%%%%%%%%%%%%%%%%%%%%%%%%%%%%%%%%%%%%%%%%%%%%%%%%%%%%%%%%%%%%%%%%%%%%%%%%%%%%%%%%%%%%%%%%%%%%%%%%%%%%%%%%%%%%%%%%%

\subsection{Morrey spaces}
\begin{defi}(Morrey spaces)
 Let us consider $q\in [1,\infty)$ and $\lambda\in [0,d)$. The Morrey space ${\cal M}_{q,\lambda}(\mathbb{R}^{d})$ is the space of functions $f\in L^{q}_{loc}(\mathbb{R}^{d})$ such that
 \begin{equation*}
  \parallel f\parallel_{{\cal M}_{q,\lambda}(\mathbb{R}^{d})}=\displaystyle{\sup_{x\in\mathbb{R}^{d}, R>0}}\,R^{-\frac{\lambda}{q}}\parallel f\cdot\chi_{B(x,R)}\parallel_{L^{q}(\mathbb{R}^{d})}\,\,\mbox{is finite}.
 \end{equation*}
 Here, $B(x,R)$ denotes an open ball in $\mathbb{R}^{d}$, centered in $x$ with radius $R$, and $\chi_{B(x,R)}$ denotes the indicator function on the ball $B(x,R)$.
 \end{defi}
The main properties on Morrey spaces are given in the next lemma. For further properties on these spaces we refer to the works \cite{Adams-Morrey-2015}, \cite{Zorko-Morrey-1986} or \cite{Sawano-Morrey-2020}.
\begin{lema}
The following itens are valid:
\begin{itemize}
\item [i.)] Let $1\leq p,q,r<\infty$ and $0\leq \lambda_{1}, \lambda_{2}, \lambda_{3}< d$ be such that $\frac{1}{r}=\frac{1}{p}+\frac{1}{q}$ and $\frac{\lambda_{3}}{r}=\frac{\lambda_{1}}{p}+\frac{\lambda_{2}}{q}$. Then, for $f\in {\cal M}_{p,\lambda_{1}}(\mathbb{R}^{d})$ and $g\in {\cal M}_{q,\lambda_{2}}(\mathbb{R}^{d})$ we have $f\cdot g \in {\cal M}_{r,\lambda_{3}}(\mathbb{R}^{d})$ and 
\begin{equation}
 \parallel f\cdot g\parallel_{{\cal M}_{r,\lambda_{3}}(\mathbb{R}^{d})}\leq \parallel f\parallel_{{\cal M}_{p,\lambda_{1}}(\mathbb{R}^{d})} \parallel g\parallel_{{\cal M}_{q,\lambda_{2}}(\mathbb{R}^{d})}.   
\end{equation}
\item[ii.)] Let $q,r\in[1,\infty)$ and  $\lambda, \mu\in [0,d)$ such that $1\leq q\leq r<\infty$ and $\frac{\lambda}{q}\geq \frac{\mu}{r}$. If $\frac{d-\lambda}{q}=\frac{d-\mu}{r}$ then it is valid the continuous inclusion
\begin{equation}
 {\cal M}_{r,\mu}(\mathbb{R}^{d})\hookrightarrow {\cal M}_{q,\lambda}(\mathbb{R}^{d}).   
\end{equation}
\item[iii.)] Let $1\leq q<\infty$ and $\lambda\in[0,d)$. If $\varphi\in L^{1}(\mathbb{R}^{d})$ and $f\in {\cal M}_{q,\lambda}(\mathbb{R}^{d})$ then $\varphi\ast f\in {\cal M}_{q,\lambda}(\mathbb{R}^{d})$ and 
\begin{equation}
 \parallel \varphi\ast f\parallel_{{\cal M}_{q,\lambda}(\mathbb{R}^{d})}\leq \parallel \varphi\parallel_{L^{1}(\mathbb{R}^{d})}\parallel f\parallel_{{\cal M}_{q,\lambda}(\mathbb{R}^{d})}.
\end{equation}
\end{itemize}   
\end{lema}
\begin{obs}
 Let us consider $q\in [1,\infty)$ and $\lambda\in [0,d)$. For each $z\in\mathbb{R}^{d}$ and $0<L$, we have
\begin{equation*}
 \parallel \chi_{B(z,L)}\parallel_{{\cal M}_{q,\lambda}(\mathbb{R}^{d})}\thickapprox_{q,d} L^{\frac{d-\lambda}{q}}.   
\end{equation*}
Here and from now on, we denote by $a\thickapprox_{q,d}b$ to say that there is a constant $C$, depending on $q$ and $d$, such that $a=Cb$.
\end{obs}

\begin{defi}
For $\vec{\lambda}\in [0,1)^{d}$, we define the space of measures of Morrey type ${\cal M}_{\vec{\lambda}}(\mathbb{R}^{d})$ as the set of Radon measures $\Upsilon$ on $\mathbb{R}^{d}$ such that
\begin{equation}
 \parallel \Upsilon\parallel_{{\cal M}_{\vec{\lambda}}(\mathbb{R}^{d})}=\displaystyle{\sup_{x\in\mathbb{R}^{d},R>0}}\,R^{-\sum_{i=1}^{d}\left(1-\lambda_{i}\right)}\mid \Upsilon\mid\left(B(x,R)\right)\,\,\,\mbox{is finite},   
\end{equation}
where $\mid \cdot \mid$ denotes the total variation of the measure $\Upsilon$.
\end{defi}

Observe that in the notation for Morrey spaces from \cite{Kozono-Yamazaki-semilinear-1994}, for $\vec{\lambda}\in[0,1)^{d}$, we have ${\cal M}^{p}(\mathbb{R}^{d})={\cal M}_{\vec{\lambda}}(\mathbb{R}^{d})$, where $p=\frac{d}{\sum_{i=1}^{d}(1-\lambda_{i})}$, and the space of Radon measures on $\mathbb{R}^{d}$ with finite total variation ${\cal M}^{1}$ is identified with ${\cal M}_{\vec{0}}$. We also observe that ${\cal M}_{\vec{q},\vec{\lambda}}(\mathbb{R}^{d})=L^{\vec{q}}(\mathbb{R}^{d})$, for $\vec{q}\in[1,\infty)^{d}$ and $\vec{\lambda}=\vec{0}$.

The next lemma is a rewritten of Lemma 1.8 from \cite{Kozono-Yamazaki-semilinear-1994}, in the language of this last definition.
\begin{lema}
Let $\nu$ denotes a Radon measure on $\mathbb{R}^{d}$ with finite total variation, $\mid \nu\mid<\infty$. Then, the following itens are valid:
\begin{itemize}
\item [i.)] For $\Upsilon\in {\cal M}_{\vec{\lambda}}(\mathbb{R}^{d})$, where $\vec{\lambda}\in [0,1)^{d}$, we have $\nu\ast\Upsilon\in {\cal M}_{\vec{\lambda}}(\mathbb{R}^{d})$ and 
\begin{equation}
\parallel \nu\ast\Upsilon\parallel_{{\cal M}_{\vec{\lambda}}(\mathbb{R}^{d})}\leq \mid \nu\mid\cdot \parallel \Upsilon\parallel_{{\cal M}_{\vec{\lambda}}(\mathbb{R}^{d})}
\end{equation}
\item [(ii.)] Let $\vec{1}\leq \vec{q}<\vec{\infty}$ and $\vec{\lambda}\in[0,1)^{d}$. If $u\in {\cal M}_{\vec{q},\vec{\lambda}}(\mathbb{R}^{d})$ then
\begin{equation}
\parallel \nu\ast u\parallel_{{\cal M}_{\vec{q}, \vec{\lambda}}(\mathbb{R}^{d})}\leq \mid \nu\mid\cdot \parallel u\parallel_{{\cal M}_{\vec{q}, \vec{\lambda}}(\mathbb{R}^{d})}
\end{equation}
\item [(iii.)] Let us assume that $\nu$ is absolutely continuous with respect to the Lebesgue measure. Then, for the measure $\Upsilon\in {\cal M}_{\vec{\lambda}}(\mathbb{R}^{d})$, the measure $\nu\ast \Upsilon$ is also absolutely continuous with respect to the Lebesgue measure and thus $\nu\ast \Upsilon\in {\cal M}_{\vec{1}, \vec{\lambda}}(\mathbb{R}^{d})$.
\end{itemize}
\end{lema}
The Besov-Morrey spaces were introduced in \cite{Kozono-Yamazaki-semilinear-1994} to study the Navier-Stokes system and the semilinear heat equations.
\begin{defi}(Besov-Morrey spaces)
\label{DefBesovType} 
Let us consider four parameters 
$$
    \text{$q\in [1,\infty)$, $\lambda\in[0,d)$, $r\in [1,\infty]$ and $\sigma\in\mathbb{R}$.}
$$    
The $d$-dimensional homogeneous Besov-Morrey spaces are defined as
\begin{itemize}
\item [(i.)] For $q\in [1,\infty)$, $\lambda\in [0,d)$, $r\in[1,\infty)$ and $\sigma\in\mathbb{R}$,
\begin{equation}
\dot{{\cal N}}^{\sigma}_{q,\lambda,r}(\mathbb{R}^{d})=\lbrace
g\in {\cal S}'/{\cal P};\parallel g\parallel_{\dot{{\cal N}}^{\sigma}_{q,\lambda,r}(\mathbb{R}^{d})}=
\left(\displaystyle{\sum_{l\in\mathbb{Z}}2^{l\sigma r}\parallel \dot{\Delta}_{l}g\parallel^{r}_{{\cal M}_{q,\lambda}(\mathbb{R}^{d})}}
\right)^{\frac{1}{r}}\,\,\mbox{is finite}
\rbrace.
\end{equation}
\item [(ii.)]  For $q\in [1,\infty)$, $\lambda\in [0,d)$, $r=\infty$ and $\sigma\in\mathbb{R}$,
\begin{equation}
\dot{{\cal N}}^{\sigma}_{q, \lambda, \infty}(\mathbb{R}^{d})=\lbrace
g\in {\cal S}'/{\cal P};\parallel g\parallel_{\dot{{\cal N}}^{\sigma}_{q,\lambda,\infty}(\mathbb{R}^{d})}=\displaystyle{\sup_{l\in\mathbb{Z}} 2^{l\sigma}\parallel \dot{\Delta}_{l}g\parallel_{{\cal M}_{q,\lambda}(\mathbb{R}^{d})}}\,\,\mbox{is finite}
\rbrace.
\end{equation}
\end{itemize}
\end{defi}%%%%%%%%%%%%%%%%%%%%%%%%%%%%%%%%%%%%%%%%%%%%%%%%%%%%%%%%%%%%%%%%%%%%%%%55
Most of the next properties on Besov-Morrey spaces were proved in \cite{Kozono-Yamazaki-semilinear-1994}.

\begin{lema}{(Main properties related to Besov-Morrey spaces)}

The following itens are valid:
\begin{itemize}
\item [(i.)] 
Let $q\in [1,\infty)$, $\lambda\in [0,d)$ and $l\in\mathbb{Z}$. There exist a constant $C$, independent on $l$, such that
\begin{equation}
 \parallel \dot{\Delta}_{l}f\parallel_{L^{\infty}(\mathbb{R}^{d})}\leq C2^{l\left(\frac{d-\lambda}{q}\right)} \parallel \dot{\Delta}_{l}f\parallel_{{\cal M}_{q,\lambda}(\mathbb{R}^{d})},
\end{equation}
for all $f\in {\cal S}'/{\cal P}$ such that $\dot{\Delta}_{l}f\in {\cal M}_{q,\lambda}(\mathbb{R}^{d})$.
\item[(ii.)] Let $q,r\in [1,\infty)$, $\lambda,\mu\in [0,d)$ and $a\in [1,\infty]$.
\begin{itemize}
    \item [(ii.1)] For $1\leq q\leq r<\infty$ and $\frac{\lambda}{q}\geq \frac{\mu}{r}$. If $\frac{d-\lambda}{q}=\frac{d-\mu}{r}$ then it is valid the continuous inclusion
    \begin{equation}
     \dot{\cal N}^{\sigma}_{r,\mu,a}(\mathbb{R}^{d})\hookrightarrow \dot{\cal N}^{\sigma}_{q,\lambda,a}(\mathbb{R}^{d}).
    \end{equation}
    \item[(ii.2)] For $b\in[1,\infty]$, if $a\leq b$ then we have the continuous inclusion
    \begin{equation}
     \dot{\cal N}^{\sigma}_{q,\lambda,a}(\mathbb{R}^{d})\hookrightarrow \dot{\cal N}^{\sigma}_{q,\lambda,b}(\mathbb{R}^{d}).
     \end{equation}
\end{itemize}
\item [(iii.)] For $q\in [1,\infty)$, $\lambda\in [0,d)$, $r\in[1,\infty]$ and $s\in\mathbb{R}$ we have the continuous inclusion
\begin{equation}
     \dot{\cal N}^{s}_{q,\lambda,r}(\mathbb{R}^{d})\hookrightarrow \dot{\cal B}^{s-\frac{d-\lambda}{q}}_{\infty,r}(\mathbb{R}^{d}).
\end{equation}
\item[(iv.)] For $q\in [1,\infty)$, $\lambda\in [0,d)$, $r\in[1,\infty]$, $s\in\mathbb{R}$ and $\theta\in(0,1)$ we have the continuous inclusion
    \begin{equation}
     \dot{\cal N}^{s}_{q,\lambda,r}(\mathbb{R}^{d})\hookrightarrow \dot{\cal N}^{s-(\frac{d-\lambda}{q})(1-\theta)}_{q/\theta,\lambda,r}(\mathbb{R}^{d}).
     \end{equation}
\item[(v.)]  For $q\in [1,\infty)$, $\lambda\in [0,d)$, $r\in[1,\infty]$ and $s\in\mathbb{R}$, the space $\dot{\cal N}^{s}_{q,\lambda,r}(\mathbb{R}^{d})$ is a Banach space and we have the continuous inclusion
\begin{equation}
 \dot{\cal N}^{s}_{q,\lambda,r}(\mathbb{R}^{d})\hookrightarrow {\cal S}'/{\cal P}.   
\end{equation}
\end{itemize}   
\end{lema}
%%%%%%%%%%%%%%%%%%%%%%%%%%%%%%%%%%%%%%%%%
The Besov type spaces we shall introduce in this work are based on the next mixed-Morrey spaces.
\subsection{Mixed-Morrey spaces}
\begin{defi}(Mixed-Morrey spaces)

 Let us consider $\vec{q}\in [1,\infty)^{d}$ and $\vec{\lambda}\in [0,1)^{d}$. The Morrey space ${\cal M}_{\vec{q},\vec{\lambda}}(\mathbb{R}^{d})$ is the space of functions $f\in L^{\vec{q}}_{\mbox{loc}}(\mathbb{R}^{d})$ such that
 \begin{equation*}
  \parallel f\parallel_{{\cal M}^{\vec{q},\vec{\lambda}}(\mathbb{R}^{d})}=\displaystyle{\sup_{x\in\mathbb{R}^{d}, R>0}}\,R^{-\left(\sum_{i=1}^{d}\frac{\lambda_{i}}{q_{i}}\right)}\parallel f\cdot\chi_{B(x,R)}\parallel_{L^{\vec{q}}(\mathbb{R}^{d})}\,\,\mbox{is finite}.
 \end{equation*}
 \end{defi}
The main properties on mixed-Morrey spaces are given in the next lemma. For further information see the works \cite{Nogayama-mixed-Morrey-2019} and \cite{Ragusa-Scapellato-mixed-Morrey-2017}.
\begin{lema}\label{properties-mixed-morrey}
The following itens are valid:
\begin{itemize}
\item [i.)] Let $\vec{p},\vec{q},\vec{r}\in [1,\infty)^{d}$ and $\vec{\lambda}_{1}, \vec{\lambda}_{2}, \vec{\lambda}_{3}\in [0,1)^{d}$ be such that $\frac{1}{r_{i}}=\frac{1}{p_{i}}+\frac{1}{q_{i}}$ and $\frac{\lambda_{3i}}{r_{i}}=\frac{\lambda_{1i}}{p_{i}}+\frac{\lambda_{2i}}{q_{i}}$, for $1\leq i\leq d$ . Then, for $f\in {\cal M}_{\vec{p},\vec{\lambda}_{1}}(\mathbb{R}^{d})$ and $g\in {\cal M}_{\vec{q},\vec{\lambda}_{2}}(\mathbb{R}^{d})$ we have $f\cdot g \in {\cal M}_{\vec{r},\vec{\lambda}_{3}}(\mathbb{R}^{d})$ and 
\begin{equation}
 \parallel f\cdot g\parallel_{{\cal M}_{\vec{r},\vec{\lambda}_{3}}(\mathbb{R}^{d})}\leq \parallel f\parallel_{{\cal M}_{\vec{p},\vec{\lambda}_{1}}(\mathbb{R}^{d})} \parallel g\parallel_{{\cal M}_{\vec{q},\vec{\lambda}_{2}}(\mathbb{R}^{d})}.   
\end{equation}
\item[ii.)] Let $\vec{q},\vec{r}\in [1,\infty)^{d}$ and  $\vec{\lambda}, \vec{\mu}\in [0,1)^{d}$ such that $\vec{1}\leq \vec{q}\leq \vec{r}<\vec{\infty}$ and $\frac{\lambda_{i}}{q_{i}}\geq \frac{\mu_{i}}{r_{i}}$ for $1\leq i\leq d$. If $\frac{1-\lambda_{i}}{q_{i}}=\frac{1-\mu_{i}}{r_{i}}$, for $1\leq i\leq d$, then it is valid the continuous inclusion
\begin{equation}
 {\cal M}_{\vec{r},\vec{\mu}}(\mathbb{R}^{d})\hookrightarrow {\cal M}_{\vec{q},\vec{\lambda}}(\mathbb{R}^{d}).   
\end{equation}
\item[iii.)] Let $\vec{1}\leq \vec{q}<\vec{\infty}$ and $\vec{\lambda}\in[0,1)^{d}$. If $\varphi\in L^{1}(\mathbb{R}^{d})$ and $f\in {\cal M}_{\vec{q},\vec{\lambda}}(\mathbb{R}^{d})$ then $\varphi\ast f\in {\cal M}_{\vec{q},\vec{\lambda}}(\mathbb{R}^{d})$ and 
\begin{equation}
 \parallel \varphi\ast f\parallel_{{\cal M}_{\vec{q},\vec{\lambda}}(\mathbb{R}^{d})}\leq \parallel \varphi\parallel_{L^{1}(\mathbb{R}^{d})}\parallel f\parallel_{{\cal M}_{\vec{q},\vec{\lambda}}(\mathbb{R}^{d})}.
\end{equation}
\item[iv.)] Let us denote by $x=(x',x'')\in \mathbb{R}^{d_{2}}$, for $x'\in\mathbb{R}^{d_{1}}$ and $x''\in\mathbb{R}^{d_{2}-d_{1}}$, where $d_{1}<d_{2}$ are positive integers. Let us consider $\vec{1}\leq \vec{q}=(q_{1},\dots, q_{d_{1}},\dots, q_{d_{2}})<+\vec{\infty}$, $\vec{0}\leq \vec{\lambda}=(\lambda_{1},\dots,\lambda_{d_{1}})<\vec{1}$ and $\vec{0}\leq \vec{\beta}=(\beta_{1},\dots,\beta_{d_{2}})<\vec{1}$. Then, every function $u(x')$ in ${\cal M}_{\vec{q},\vec{\lambda}}(\mathbb{R}^{d_{1}})$ is also an element in ${\cal M}_{\vec{q},\vec{\beta}}(\mathbb{R}^{d_{2}})$, if the next condition
\begin{equation}\label{condition-Morrey}
 \displaystyle{\sum_{i=1}^{d_{2}}}\frac{1-\beta_{i}}{q_{i}}=\displaystyle{\sum_{i=1}^{d_{1}}}\frac{1-\lambda_{i}}{q_{i}}   
\end{equation}
is valid.
\end{itemize}

\end{lema}
\begin{obs}
Let us have in mind the next itens:
\begin{itemize}
\item The item iv.) from the last Lemma \ref{properties-mixed-morrey} is new and its proof follows the computations from (\cite{Kozono-Yamazaki-semilinear-1994}, Proposition 1.3). 
\item Let us consider $\vec{q}\in [1,\infty)^{d}$ and $\vec{\lambda}\in [0,1)^{d}$. For each $z\in\mathbb{R}^{d}$ and $0<L$ we have
\begin{equation*}
 \parallel \chi_{B(z,L)}\parallel_{{\cal M}_{\vec{q},\vec{\lambda}}} \thickapprox_{\vec{q},d}   L^{\sum_{i=1}^{d}\left(\frac{1-\lambda_{i}}{q_{i}}\right)}. 
\end{equation*}
\end{itemize}
\end{obs}
\subsection{Fourier-Besov mixed-Lebesgue spaces}
Fourier-Besov spaces were introduced in the work of P. Konieczny and T. Yoneda \cite{Konieczny-Yoneda-2011} to study the dispersive effect of the Coriolis force for both the stationary and the non-stationary Navier-Stokes system. 

Next spaces were employed in \cite{Leithold-Neves-Besov-mixed-Lebesgue} in order to study the Navier-Stokes system.
%%%%%%%%%%%%%%%%%%%%%%%%%%%%%%%%%%%%%%%%%%%%

\begin{defi}(Fourier-Besov mixed-Lebesgue spaces)\label{DefFourier-BesovType}
Let us consider the three parameters 
$$
    \text{$\vec{p}=(p_1,\dots,p_d)\in [1,\infty]^{d}$, $q\in [1,\infty]$ and $s\in\mathbb{R}$.}
$$ 
The $d$-dimensional homogeneous Fourier-Besov mixed-Lebesgue spaces are defined as
\begin{itemize}
\item [(i.)] For $\vec{p}\in [1,\infty]^{d}$, $q\in[1,\infty)$ and $s\in\mathbb{R}$,
\begin{equation}
\dot{{\cal FB}}^{s}_{\vec{p},q}(\mathbb{R}^{d})=\lbrace
g\in {\cal S}'/{\cal P};\parallel g\parallel_{\dot{{\cal FB}}^{s}_{\vec{p},q}(\mathbb{R}^{d})}=
\left(\displaystyle{\sum_{l\in\mathbb{Z}}2^{lsq}\parallel \phi_{l}\hat{g}\parallel^{q}_{L^{\vec{p}}(\mathbb{R}^{d})}}
\right)^{\frac{1}{q}}\,\,\mbox{is finite}
\rbrace.
\end{equation}
\item [(ii.)]  For $\vec{p}\in [1,\infty]^{d}$, $q=\infty$ and $s\in\mathbb{R}$,
\begin{equation}
\dot{{\cal FB}}^{s}_{\vec{p},\infty}(\mathbb{R}^{d})=\lbrace
g\in {\cal S}'/{\cal P};\parallel g\parallel_{\dot{{\cal FB}}^{s}_{\vec{p},\infty}(\mathbb{R}^{d})}=\displaystyle{\sup_{l\in\mathbb{Z}} 2^{ls}\parallel \phi_{l}\hat{g}\parallel_{L^{\vec{p}}(\mathbb{R}^{d})}}\,\,\mbox{is finite}
\rbrace. 
\end{equation}
\end{itemize}
We also called these spaces as the $d$-dimensional homogeneous mixed-norm Fourier-Besov-Lebesgue spaces or the $d$-dimensional homogeneous Fourier-Besov spaces based on mixed-Lebesgue spaces.
\end{defi}

The next properties were established in our previous work $\cite{Leithold-Neves-Besov-mixed-Lebesgue}$.
\begin{lema}{(Main properties related to Fourier-Besov mixed-Lebesgue spaces)}\label{properties-Fourier-Besov-mixed}

The following itens are valid.
\begin{itemize}
\item [(i.)] Let $\vec{p},\vec{q}\in[1,\infty]^{d}$ such that $\vec{p}\leq \vec{q}$ and let $g$ be any map with $\hat{g}\in L^{\vec{q}}(\mathbb{R}^{d})$ such that 

\begin{equation*}
 \supp \hat{g}\subset\{
 \xi_{1}\in\mathbb{R};\mid \xi_{1}\mid \leq A_{1}2^{i_1}\}\times\cdots\times\{
 \xi_{d}\in\mathbb{R};\mid \xi_{d}\mid \leq A_{d}2^{i_d}\},   
\end{equation*}
for some real constants $A_1,\cdots, A_d$ and integers $i_{1},\cdots, i_{d}\in \mathbb{Z}$. Then,
\begin{equation*}
 \parallel\left(\xi_1\right)^{\beta_1}\cdots \left(\xi_d\right)^{\beta_d}\hat{g}\parallel_{L^{\vec{p}}(\mathbb{R}^{d})}\leq C2^{k_{0}}  \parallel\hat{g}\parallel_{L^{\vec{q}}(\mathbb{R}^{d})},
\end{equation*}
where $k_{0}=\displaystyle{\sum_{k=1}^{d}}i_{k}\left(\mid \beta_{k}\mid+\frac{1}{p_{k}}-\frac{1}{q_{k}}\right)$ and $C$ is a constant depending on $A_{1},\cdots,A_{d}$, for any indexes $\beta_{1},\cdots,\beta_{d}$.
\item[(ii.)] Let us consider $\vec{p},\vec{q}\in[1,\infty]^{d}$, $\vec{q}\leq\vec{p}$ and $1\leq a_{1}\leq a_{2}\leq \infty$. If we consider $s_{1},s_{2}\in\mathbb{R}$ such that $s_{1}\geq s_{2}$ and $s_{1}+\displaystyle{\sum_{i=1}^{d}}\frac{1}{p_{i}}=s_{2}+\displaystyle{\sum_{i=1}^{d}}\frac{1}{q_{i}}$, then we have
\begin{equation*}
 \dot{{\cal FB}}^{s_{1}}_{\vec{p},a_{1}}(\mathbb{R}^{d})\hookrightarrow
 \dot{{\cal FB}}^{s_{2}}_{\vec{q},a_{2}}(\mathbb{R}^{d}),
\end{equation*}
that is, the above inclusion is continuous.
\item [(iii.)] For $\vec{p}\in[1,\infty]^{d}$, $q\in[1,\infty]$ and $s\in\mathbb{R}$ the Fourier-Besov mixed-Lebesgue space  $\dot{{\cal FB}}^{s}_{\vec{p},q}(\mathbb{R}^{d})$ is a Banach space endowed with the norm $\parallel \cdot\parallel_{\dot{{\cal FB}}^{s}_{\vec{p},q}(\mathbb{R}^{d})}$.
\end{itemize}   
\end{lema}
%%%%%%%%%%%%%%%%%%%%%%%%%%%%%%%%%%%%%%%%%%%%%%%%
The next spaces were employed in \cite{Almeida-Ferreira-Lima-uniform-2017} to study the Navier-Stokes-Coriolis system.
\begin{defi}(Fourier-Besov-Morrey spaces) 
Let us consider four parameters
$$
    \text{$q\in[1,\infty)$, $\lambda\in[0,d)$, $r\in[1,\infty]$ and $s\in\mathbb{R}$.}
$$ 
The $d$-dimensional homogeneous Fourier-Besov-Morrey spaces are defined as
\begin{itemize}
\item [(i.)] For $q\in[1,\infty)$, $\lambda\in[0,d)$, $r\in[1,\infty)$ and $s\in\mathbb{R}$,
\begin{equation}
\dot{{\cal FN}}^{s}_{q,\lambda,r}(\mathbb{R}^{d})=\lbrace
g\in {\cal S}'/{\cal P};\parallel g\parallel_{\dot{{\cal FN}}^{s}_{q,\lambda,r}(\mathbb{R}^{d})}=
\left(\displaystyle{\sum_{l\in\mathbb{Z}}2^{lsr}\parallel \phi_{l}\hat{g}\parallel^{r}_{{\cal M}_{q,\lambda}(\mathbb{R}^{d})}}
\right)^{\frac{1}{r}}\,\,\mbox{is finite}
\rbrace.
\end{equation}
\item [(ii.)]  For  $q\in[1,\infty)$, $\lambda\in[0,d)$, $r=\infty$ and $s\in\mathbb{R}$,
\begin{equation}
\dot{{\cal FN}}^{s}_{q,\lambda,\infty}(\mathbb{R}^{d})=\lbrace
g\in {\cal S}'/{\cal P};\parallel g\parallel_{\dot{{\cal FN}}^{s}_{q,\lambda,\infty}(\mathbb{R}^{d})}=\displaystyle{\sup_{l\in\mathbb{Z}} 2^{ls}\parallel \phi_{l}\hat{g}\parallel_{{\cal M}_{q,\lambda}(\mathbb{R}^{d})}}\,\,\mbox{is finite}
\rbrace. 
\end{equation}
\end{itemize}
We also called these spaces as the $d$-dimensional homogeneous mixed-norm Besov-Morrey spaces or the $d$-dimensional homogeneous Fourier-Besov spaces based on Morrey spaces.
\end{defi}
%%%%%%%%%%%%%%%%%%%%%%%%%%%%%%%%%%%%%%%%%%%%%%%%
Most of the next properties were proved in \cite{Almeida-Ferreira-Lima-uniform-2017}, (see also the references therein). 
\begin{lema}{(Main properties related to Fourier-Besov-Morrey spaces)}

\begin{itemize}
\item [(i.)] 
Let $q\in [1,\infty)$, $\lambda\in [0,d)$ and $l\in\mathbb{Z}$. There exist a constant $C$, independent on $l$, such that
\begin{equation}
 \parallel \phi_{l} \hat{f}\parallel_{L^{1}(\mathbb{R}^{d})}\leq C\cdot2^{l\left(d-\frac{d-\lambda}{q}\right)} \parallel \phi_{l} \hat{f}\parallel_{{\cal M}_{q,\lambda}(\mathbb{R}^{d})},
\end{equation}
for all $f\in {\cal S}'/{\cal P}$ such that $\phi_{l} \hat{f}\in {\cal M}_{q,\lambda}(\mathbb{R}^{d})$.
\item[(ii.) ] 
Let $1\leq q<r<\infty$ and $\lambda,\mu\in[0,d)$ be parameters such that $\frac{\mu}{r}\leq \frac{\lambda}{q}$ and $\frac{d-\mu}{r}< \frac{d-\lambda}{q}$, and $l\in\mathbb{Z}$. Then, there exist a positive constant $C$, independent on $l$, such that
\begin{equation}
 \parallel \phi_{l}\hat{f}\parallel_{{\cal M}_{q,\lambda}(\mathbb{R}^{d})}\leq C\cdot2^{l\left(\frac{d-\lambda}{q}-\frac{d-\mu}{r}\right)}\parallel \phi_{l}\hat{f}\parallel_{{\cal M}_{r,\mu}(\mathbb{R}^{d})},   
\end{equation}
for each $f\in{\cal S}'/{\cal P}$ such that $\phi_{l}\hat{f}\in{\cal M}_{r,\mu}(\mathbb{R}^{d})$.
\item[(iii.)] Let $q,r\in [1,\infty)$, $\lambda,\mu\in [0,d)$ and $a\in [1,\infty]$.
\begin{itemize}
    \item [(iii.1)] For $1\leq q\leq r<\infty$ and $\frac{\lambda}{q}\geq \frac{\mu}{r}$. If $\frac{d-\lambda}{q}=\frac{d-\mu}{r}$ then it is valid the continuous inclusion
    \begin{equation}
     \dot{\cal FN}^{\sigma}_{r,\mu,a}(\mathbb{R}^{d})\hookrightarrow \dot{\cal FN}^{\sigma}_{q,\lambda,a}(\mathbb{R}^{d}).
    \end{equation}
    \item[(iii.2)] For $b\in[1,\infty]$, if $a\leq b$ then we have the continuous inclusion
    \begin{equation}
     \dot{\cal FN}^{\sigma}_{q,\lambda,a}(\mathbb{R}^{d})\hookrightarrow \dot{\cal FN}^{\sigma}_{q,\lambda,b}(\mathbb{R}^{d}).
     \end{equation}
\end{itemize}
\item[(iv.)] Let $1\leq q<r<\infty$, $\lambda,\mu\in[0,d)$, $a\in[1,\infty]$ and $s_{1},s_{2}\in\mathbb{R}$ be parameters such that $\frac{\mu}{r}\leq \frac{\lambda}{q}$ and $\frac{d-\mu}{r}< \frac{d-\lambda}{q}$, and $s_{2}<s_{1}$. If $s_{2}+\frac{d-\lambda}{q}=s_{1}+\frac{d-\mu}{r}$ then
\begin{equation}
\dot{\cal FN}^{s_{1}}_{r,\mu,a}(\mathbb{R}^{d})\hookrightarrow \dot{\cal FN}^{s_{2}}_{q,\lambda,a}(\mathbb{R}^{d}).
\end{equation}
\item [(v.)] For $q\in [1,\infty)$, $\lambda\in [0,d)$, $r\in[1,\infty]$ and $\sigma\in\mathbb{R}$ we have the next continuous inclusions
\begin{itemize}
\item [(v.1)]
\begin{equation}
\dot{\cal FN}^{\sigma}_{q,\lambda,r}(\mathbb{R}^{d})\hookrightarrow \dot{\cal FB}^{\sigma-\left(d-\frac{d-\lambda}{q}\right)}_{1,r}(\mathbb{R}^{d}),\,\,\mbox{and}
\end{equation}
\item [(v.2)]
\begin{equation}
\dot{{\cal FB}}^{\sigma}_{1,r}(\mathbb{R}^{d})\hookrightarrow \dot{{\cal B}}^{\sigma}_{\infty,r}(\mathbb{R}^{d}).    
\end{equation}
\end{itemize}

\item[(vi.)] For $q\in [1,\infty)$, $\lambda\in [0,d)$, $r\in[1,\infty]$, $s\in\mathbb{R}$ and $\theta\in(0,1)$ we have the continuous inclusion
    \begin{equation}
     \dot{\cal FN}^{s}_{q,\lambda,r}(\mathbb{R}^{d})\hookrightarrow \dot{\cal FN}^{s-(\frac{d-\lambda}{q})(1-\theta)}_{q/\theta,\lambda,r}(\mathbb{R}^{d}).
     \end{equation}
\item[(vii.)]  For $q\in [1,\infty)$, $\lambda\in [0,d)$, $r\in[1,\infty]$ and $s\in\mathbb{R}$, the space $\dot{\cal FN}^{s}_{q,\lambda,r}(\mathbb{R}^{d})$ is a Banach space and we have the continuous inclusion
\begin{equation}
 \dot{\cal FN}^{s}_{q,\lambda,r}(\mathbb{R}^{d})\hookrightarrow {\cal S}'/{\cal P}.   
\end{equation}
\end{itemize}   
\end{lema}

\subsection{The Bony's paraproduct and a fixed point lemma}
%%%%%%%%%%%%%%%%%%%

In this subsection we recall two main analytical tools which we shall use along this work, as being the Bony's 
paraproduct and a fixed point scheme. 

Let us consider the Bony's paraproduct decomposition. Given tempered distributions $u,v\in {\cal S}'/{\cal P}$, we consider the decomposition
\begin{equation*}
 vw=T_{v}w+T_{w}v+R(v,w),   
\end{equation*}
where the Bony's paraproduct operator $T_v(\cdot)$ is given by 
\begin{equation*}
T_{v}(w)=\displaystyle{\sum_{l\in\mathbb{Z}}}\dot{S}_{l-1}v \dot{\Delta}_{l}w,   
\end{equation*}
and the remainder part $R(\cdot,\cdot)$ is given by 
\begin{equation*}
R(v,w)=\displaystyle{\sum_{l\in\mathbb{Z}}}\dot{\Delta}_{l}v\tilde{\dot{\Delta}}_{l}w,\,\,\mbox{where},\,\,
\tilde{\dot{\Delta}}_{l}w=\displaystyle{\sum_{\mid l-l'\mid\leq 1}}\dot{\Delta}_{l'}w.
\end{equation*}
For further information on 
Bony's paraproduct see (\cite{Lemarie-NS-century-2018}, Chapter 3), (\cite{Bahouri-Chemin-Danchin-Fourier-PDE-2011}, 
Section 2.6) and the original work by J. M. Bony \cite{Bony-symbolique-1981}.

%%%%%%%%%%%%%%%%%%
%%%%%%%%%%%%%%%%%%

Now we recall an abstract lemma related to a fixed point argument. First, we note that the mild formulation has the abstract form 
\begin{equation}\label{mild-formulation-abstract}
u(t)=z_{0}+B(u,u)(t),     
\end{equation}
for all $t\in [0,\infty)$, where $z_{0}=S_{\nu}(t)u_{0}$ and $B(\cdot,\cdot)$ is the given bilinear operator. 

The following lemma (\cite{Lemarie-NS-century-2018}, Theorem 13.2) is the fixed point scheme which we shall apply in order to solve the abstract equation associated to the mild formulation for the Navier-Stokes system.
\begin{lema}\label{fixed-point-scheme}
Let $\left(Z,\parallel \cdot\parallel\right)$ be a Banach space and $B:Z\times Z\rightarrow Z$ be a bilinear operator for which there is $K>0$ such that $\parallel B\parallel_{{\cal B}(Z)}\leq K$. Then, for $0<\varepsilon<\frac{1}{4K}$ and $\parallel z_{0}\parallel_{Z}\leq \varepsilon$, the abstract equation $z=z_{0}+B(z,z)$ has a unique solution in the closed ball $B_{2\varepsilon}(0)=\{z\in Z;\parallel z\parallel_{Z}\leq 2\varepsilon\}$ and $\parallel z\parallel_{Z}\leq 2\parallel z_{0}\parallel_{Z}$. Moreover, the solution depends continuously on initial data, in the following sense: for $\tilde{z}_{0}\in Z$, with $\parallel\tilde{z}_{0}\parallel_{Z}\leq \varepsilon$, if $\tilde{z}=\tilde{z}_{0}+B(\tilde{z},\tilde{z})$ and $\parallel \tilde{z}\parallel_{Z}\leq 2\varepsilon$, then 
\begin{equation*}
 \parallel z-\tilde{z}\parallel_{Z}\leq \left(1-4K\varepsilon\right)^{-1}\parallel z_{0}-\tilde{z}_{0}\parallel_{Z}.   
\end{equation*}
\end{lema}

\section{Besov mixed-Morrey spaces}\label{section-BMMorrey}
Based on the above mixed-Morrey spaces, now we introduce the Besov mixed-Morrey spaces. These new spaces are important to study the Navier-Stokes system.
\begin{defi}(Besov mixed-Morrey spaces)
\label{DefBesovType} 

Let us consider four parameters 
$$
    \text{$\vec{q}\in [1,\infty)^{d}$, $\vec{\lambda}\in[0,1)^{d}$, $r\in [1,\infty]$ and $\sigma\in\mathbb{R}$.}
$$    
The $d$-dimensional homogeneous Besov mixed-Morrey spaces are defined as
\begin{itemize}
\item [(i.)] For $\vec{q}\in [1,\infty)^{d}$, $\vec{\lambda}\in [0,1)^{d}$, $r\in[1,\infty)$ and $\sigma\in\mathbb{R}$,
\begin{equation}
\dot{{\cal N}}^{\sigma}_{\vec{q},\vec{\lambda},r}(\mathbb{R}^{d})=\lbrace
g\in {\cal S}'/{\cal P};\parallel g\parallel_{\dot{{\cal N}}^{\sigma}_{\vec{q},\vec{\lambda},r}(\mathbb{R}^{d})}=
\left(\displaystyle{\sum_{l\in\mathbb{Z}}2^{l\sigma r}\parallel \dot{\Delta}_{l}g\parallel^{r}_{{\cal M}_{\vec{q},\vec{\lambda}}(\mathbb{R}^{d})}}
\right)^{\frac{1}{r}}\,\,\mbox{is finite}
\rbrace.
\end{equation}
\item [(ii.)]  For $q\in [1,\infty)$, $r=\infty$ and $\sigma\in\mathbb{R}$,
\begin{equation}
\dot{{\cal N}}^{\sigma}_{\vec{q}, \vec{\lambda},\infty}(\mathbb{R}^{d})=\lbrace
g\in {\cal S}'/{\cal P};\parallel g\parallel_{\dot{{\cal N}}^{\sigma}_{\vec{q},\vec{\lambda},\infty}(\mathbb{R}^{d})}=\displaystyle{\sup_{l\in\mathbb{Z}} 2^{l\sigma}\parallel \dot{\Delta}_{l}g\parallel_{{\cal M}_{\vec{q},\vec{\lambda}}(\mathbb{R}^{d})}}\,\,\mbox{is finite}
\rbrace.
\end{equation}
\end{itemize}
We also called these spaces as the $d$-dimensional homogeneous mixed-norm Besov-Morrey spaces or the $d$-dimensional homogeneous Besov spaces based on mixed-Morrey spaces.
\end{defi}
%%%%%%%%%%%%%%%%%%%%%%%%%%%%%%%%%%%%%%%%%%%%5
\subsection{Proof of Theorem \ref{properties-Besov-spaces-mixed-Morrey}}
\begin{proof}
\begin{itemize}
\item [(i.)] First, we denote $K_{0}={\cal F}^{-1}\tilde{\phi}_{0}$, where $\tilde{\phi}_{0}=\phi_{-1}+\phi_{0}+\phi_{1}$, we recall that $\tilde{\phi}_{l}=\phi_{l-1}+\phi_{l}+\phi_{l+1}$ and we observe that ${\cal F}^{-1}\left(\tilde{\phi}_{l}\right)=2^{ld}\left( {\cal F}^{-1}\phi_{0}\right)(2^{l}\cdot)$. Then we have the estimates
\begin{equation}\label{split0}
\begin{split}
 \mid \dot{\Delta}_{l}u(x)\mid &=\mid{\cal F}^{-1}[\tilde{\phi}_{l}\phi_{l}\hat{u}](x)\mid \\
    &=\mid{\cal F}^{-1}(\tilde{\phi}_{l})\ast {\cal F}^{-1}(\phi_{l}\hat{u})(x)\mid\\
    &=\mid \left(K_{l}\ast \dot{\Delta}_{l}u\right)(x)\mid\\
    &=\mid \displaystyle{\int_{\mathbb{R}^{d}}}K_{l}(y)\dot{\Delta}_{l}u(x-y)\,dy\mid\\
    &\leq A_{1}+A_{2}
\end{split}    
\end{equation}
where
\begin{equation*}
A_{1}= \displaystyle{\int_{\mid y\mid\leq 2^{-l}}}\mid K_{l}(y)\mid \mid\dot{\Delta}_{l}u(x-y)\mid\,dy   
\end{equation*}
and 
\begin{equation*}
A_{2}= \displaystyle{\sum_{k=-l+1}^{\infty}}\displaystyle{\int_{2^{k-1}\leq\mid y\mid\leq 2^{-k}}}\mid K_{l}(y)\mid \mid\dot{\Delta}_{l}u(x-y)\mid\,dy.   
\end{equation*}
 The estimate for the term $A_{1}$ is given by rewritten the integral by mean of a suitable indicator function $\chi_{B}$ and by using the H\"older's inequality for mixed-Lebesgue spaces we get
 \begin{equation*}
  \begin{split}
   A_{1}&=\displaystyle{\int_{\mathbb{R}^{d}}}\mid K_{l}(y)\mid \mid\dot{\Delta}_{l}u(x-y)\mid\chi_{\{y;\mid y\mid\leq 2^{-l}\}}(y)\,dy\\
   &\leq \parallel K_{l}\parallel_{L^{\vec{q}'}(\mathbb{R}^{d})}\parallel \dot{\Delta}_{l}u(x-\cdot)\chi_{\{y;\mid y\mid\leq 2^{-l}\}}(\cdot)\parallel_{L^{\vec{q}}(\mathbb{R}^{d})}\\
   &=\parallel K_{l}\parallel_{L^{\vec{q}'}(\mathbb{R}^{d})}\parallel \dot{\Delta}_{l}u \cdot\chi_{\{z;\mid z-x\mid\leq 2^{-l}\}}(\cdot)\parallel_{L^{\vec{q}}(\mathbb{R}^{d})}\\
   &=2^{l(d-\sum_{i=1}^{d}\frac{1}{q_{i}'})}\parallel K_{0}\parallel_{L^{\vec{q}'}(\mathbb{R}^{d})}\parallel\dot{\Delta}_{l}u\cdot \chi_{B(x,2^{-l})}(\cdot)\parallel_{L^{\vec{q}}(\mathbb{R}^{d})},
  \end{split}   
 \end{equation*}
 where we also used the change of variables $z=x-y$. Therefore, by definition of mixed-Morrey spaces ${\cal M}_{\vec{q},\vec{\lambda}}(\mathbb{R}^{d})$, we obtain
 \begin{equation}\label{equationA1}
  A_{1}\leq 2^{l(\sum_{i=1}^{d}\frac{1-\lambda_{i}}{q_{i}})}\parallel \dot{\Delta}_{l}u\parallel_{{\cal M}_{\vec{q},\vec{\lambda}}(\mathbb{R}^{d})}\parallel K_{0}\parallel_{L^{\vec{q}'}(\mathbb{R}^{d})}.
\end{equation}
Now we compute the estimate for the term $A_{2}$. We first rewrite the term $A_{2}$ and by the H\"older's inequality for mixed-Lebesgue spaces we get the next estimate
\begin{equation*}
 \begin{split}
 A_{2}&= \displaystyle{\sum_{k=-l+1}^{\infty}}\displaystyle{\int_{\mathbb{R}^{d}}}\mid K_{l}(y)\chi_{\{y;\mid y\mid\geq 2^{k-1}\}}(y)\mid \mid\dot{\Delta}_{l}u(x-y)\mid\chi_{\{y;\mid y\mid\leq 2^{k}\}}\,dy\\
 &\leq \displaystyle{\sum_{k=-l+1}^{\infty}} \parallel K_{l}(\cdot)\chi_{\{y;\mid y\mid\geq 2^{k-1}\}}(\cdot)\parallel_{L^{\vec{q}'}(\mathbb{R}^{d})}\parallel \dot{\Delta}_{l}u \cdot\chi_{B(x,2^{k})}\parallel_{L^{\vec{q}}(\mathbb{R}^{d})}\\
 &=\displaystyle{\sum_{k=-l+1}^{\infty}} 2^{l(d-\sum_{i=1}^{d}\frac{1}{q'_{i}})} \parallel K_{0}(\cdot)\chi_{\{z;\mid z\mid\geq 2^{k+l-1}\}}(\cdot)\parallel_{L^{\vec{q}'}(\mathbb{R}^{d})}\parallel \dot{\Delta}_{l}u \cdot\chi_{B(x,2^{k})}\parallel_{L^{\vec{q}}(\mathbb{R}^{d})},
 \end{split}   
\end{equation*}
where we also used the change of variables $z=2^{l}y$. Therefore, by definition of mixed-Morrey spaces ${\cal M}_{\vec{q},\vec{\lambda}}(\mathbb{R}^{d})$ and introducing the new index $m=k+l-1$, we get
\begin{equation}\label{equationA2}
 A_{2}\leq 2^{l(\sum_{i=1}^{d}\frac{1-\lambda_{i}}{q_{i}})}\cdot  \left(
 \displaystyle{\sum_{m=0}^{\infty}}2^{m(\sum_{i=1}^{d}\frac{\lambda_{i}}{q_{i}})}A_{2m}
 \right)\parallel \dot{\Delta}_{l}u\parallel_{{\cal M}_{\vec{q},\vec{\lambda}}(\mathbb{R}^{d})},
\end{equation}
where 
\begin{equation*}
  A_{2m}= \parallel K_{0}(\cdot)\chi_{B^{c}(0,2^{m})}(\cdot)\parallel _{L^{\vec{q}'}(\mathbb{R}^{d})}.
\end{equation*}
Now, we shall establish that
\begin{equation}\label{equationA2m}
  A_{2m}\leq  2^{-m(L_{1}-\frac{1}{q_{1}'})}+\dots+2^{-m(L_{d}-\frac{1}{q_{d}'})}, 
\end{equation}
where $L_{j}$ is a positive integer such that $L_{j}>\frac{1}{q_{j}'}$, for each $j=1,\dots,d$.
First, we note that 
\begin{equation*}
 A_{2m}\leq \displaystyle{\sum_{j=1}^{d}}\parallel K_{0}(\cdot)\chi_{\{x;\mid x_{j}\mid>\frac{2^{m}}{\sqrt{d}}\}\}}(\cdot)   \parallel_{L^{\vec{q}'}(\mathbb{R}^{d})}=\displaystyle{\sum_{j=1}^{d}}B_{j},
\end{equation*}
since $B^{c}(0,2^{m})\subset \mathbb{R}^{d}-[-a,a]^{d}$, for $a=\frac{2^{m}}{\sqrt{d}}$. Let us note that
\begin{equation*}
B_{j}\leq \parallel \frac{1}{(1+\mid x_{1}\mid)^{L_{1}}}\parallel_{L^{q_{1}'}}\dots\parallel \frac{1}{(1+\mid x_{j}\mid)^{L_{j}}}\chi_{\{x_{j};\mid x_{j}\mid>a\}}(x_{j})\parallel_{L^{q_{j}'}}\dots \parallel \frac{1}{(1+\mid x_{d}\mid)^{L_{d}}}\parallel_{L^{q_{d}'}}\cdot\parallel K_{0}\parallel_{L,{\cal S}},
\end{equation*}
for $L=L_{1}+ \dots + L_{d}$ integer, since each one $L_{i}$ is also a positive integer. By easy computations we get the estimate
\begin{equation*}
 B_{j}\leq C(L_{1},\dots,L_{d},q_{1}',\dots,q_{d}',d) 2^{-m(L_{j}-\frac{1}{q_{j}'})},
\end{equation*}
where each $L_{j}$ is choosen a positive integer such that $L_{j}>\frac{1}{q_{j}'}$, for each $j=1,\dots,d$, and so (\ref{equationA2m}) is proved. By considering (\ref{equationA2m}) in (\ref{equationA2}) we get 
\begin{equation}\label{series-converge}
\begin{split}
 A_{2}&\leq 2^{l(\sum_{i=1}^{d}\frac{1-\lambda_{i}}{q_{i}})}\left(\displaystyle{\sum_{j=1}^{d}} \displaystyle{\sum_{m=0}^{\infty}}2^{m(\sum_{i=1}^{d}\frac{\lambda_{i}}{q_{i}} -(L_{j}-\frac{1}{q_{j}'}))}\right)\cdot\parallel K_{0}\parallel_{L,{\cal S}}\parallel \dot{\Delta}_{l}u\parallel_{{\cal M}_{\vec{q},\vec{\lambda}}(\mathbb{R}^{d})}.
\end{split}
\end{equation}
For each $j=1,\dots, d$, if we choose an integer $L_{j}$ such that 
\begin{equation*}
  L_{j}>\frac{1}{q_{j}'}+ \sum_{i=1}^{d}\frac{\lambda_{i}}{q_{i}},
\end{equation*}
then the series in (\ref{series-converge}) converges. Therefore we get 
\begin{equation*}
   A_{2}\leq C 2^{l(\sum_{i=1}^{d}\frac{1-\lambda_{i}}{q_{i}})}\parallel K_{0}\parallel_{L,{\cal S}}\parallel \dot{\Delta}_{l}u\parallel_{{\cal M}_{\vec{q},\vec{\lambda}}(\mathbb{R}^{d})}.
\end{equation*}
If we replace this last estimate together with (\ref{equationA1}) in (\ref{split0}) we obtain the estimate
\begin{equation*}
    \mid \dot{\Delta}_{l}u(x)\mid\leq C 2^{l(\sum_{i=1}^{d}\frac{1-\lambda_{i}}{q_{i}})}\parallel \dot{\Delta}_{l}u\parallel_{{\cal M}_{\vec{q},\vec{\lambda}}(\mathbb{R}^{d})},
\end{equation*}
which implies the desired estimate.
\item [(ii.)] To prove item (ii.1) it is enough to consider property (ii.) from Lemma \ref{properties-mixed-morrey}, and to prove item (ii.2) it is enough to recall the continuous inclusion $l^{a}\hookrightarrow l^{b}$.
\item [(iii.)] In order to prove item (iii) it is enough to use the above Bernstein type inequality, item (i). 
\item [(iv.)] The proof of item (iv.) follows the lines from Theorem 2.4 from  \cite{Kozono-Yamazaki-semilinear-1994}.
\item [(v.)] The proof of item (v.) is given similarly to the proof of Corollary 2.6 from \cite{Kozono-Yamazaki-semilinear-1994}.
\item [(vi.)] The proof follows from the direct application of item iv.) from Lemma \ref{properties-mixed-morrey}.  
\end{itemize}
\end{proof}

%%%%%%%%%%%%%%%%%%%%%%%%%%%%%%%%%%%%%%%%%%%%%%%%%%%
The next lemma is an extension of Proposition 2.1 from \cite{Kozono-Yamazaki-semilinear-1994}. 
\begin{prop}
There is a positive constant $D>0$ such that the following is valid.
Given $\gamma\in\mathbb{R}$ and $l\in\mathbb{Z}$. Let $P(\eta)$ be a $C^{\infty}$-function on $\tilde{{\cal A}}_{l}={\cal A}_{l-1}\cup {\cal A}_{l}\cup {\cal A}_{l+1}$ such that it is valid the estimate
\begin{equation}
 \mid \left( \partial^{\mid \alpha\mid}/\partial \eta^{\alpha} \right)P(\eta)\mid \leq A2^{\left(\gamma-\mid \alpha\mid\right)l}, \,\,\mbox{for all}\,\eta\in \tilde{{\cal A}}_{l},   
\end{equation}
for some constant $A$, for every $\alpha\in \mathbb{N}^{d}$ such that $\mid \alpha\mid \leq [d/2]+1$. Suppose, $\vec{q}\in[1,\infty)^{d}$ and $\vec{\lambda}\in[0,1)^{d}$. Then, for each $u\in {\cal M}_{\vec{q},\vec{\lambda}}(\mathbb{R}^{d})$ with $\supp\hat{u}\subset{{\cal A}_{l}}$, we have ${\cal F}^{-1}[P(\eta)\hat{u}]\in {\cal M}_{\vec{q},\vec{\lambda}}(\mathbb{R}^{d})$ and
\begin{equation}
\parallel {\cal F}^{-1}[P(\eta)\hat{u}]\parallel_{{\cal M}_{\vec{q},\vec{\lambda}}(\mathbb{R}^{d})}\leq DA2^{\gamma l}\parallel u\parallel_{{\cal M}_{\vec{q},\vec{\lambda}}(\mathbb{R}^{d})}.
\end{equation}
\begin{proof}
The proof is similar to the proof of Proposition 2.1 from \cite{Kozono-Yamazaki-semilinear-1994}, since it is valid the Young's inequality for convolution, item (iii.) from Lemma \ref{properties-mixed-morrey}.
\end{proof}
\end{prop}
The next lemma is an extension of Proposition 2.2 from \cite{Kozono-Yamazaki-semilinear-1994}. 
\begin{prop}
There is a positive constant $D>0$ such that the following is valid. Let $P(\eta)$ is a $C^{\infty}$-function on $B(0,\frac{8}{3})$ such that it is valid the estimate
\begin{equation}
 \mid \left( \partial^{\mid \alpha\mid}/\partial \eta^{\alpha} \right)P(\eta)\mid \leq A, \,\,\mbox{for all}\,\eta\in B(0,\frac{8}{3}),   
\end{equation}
for some constant $A$, for every $\alpha\in \mathbb{N}^{d}$ such that $\mid \alpha\mid \leq [d/2]+1$. Suppose, $\vec{q}\in[1,\infty)^{d}$ and $\vec{\lambda}\in[0,1)^{d}$. Then, for each $u\in {\cal M}_{\vec{q},\vec{\lambda}}(\mathbb{R}^{d})$ with $\supp\hat{u}\subset B(0,\frac{4}{3})$, we have ${\cal F}^{-1}[P(\eta)\hat{u}]\in {\cal M}_{\vec{q},\vec{\lambda}}(\mathbb{R}^{d})$ and
\begin{equation}
\parallel {\cal F}^{-1}[P(\eta)\hat{u}]\parallel_{{\cal M}_{\vec{q},\vec{\lambda}}(\mathbb{R}^{d})}\leq DA\parallel u\parallel_{{\cal M}_{\vec{q},\vec{\lambda}}(\mathbb{R}^{d})}.
\end{equation}
\end{prop}
\begin{proof}
The proof follow the same reasoning in Proposition 2.2 from \cite{Kozono-Yamazaki-semilinear-1994}, because it is valid the Young's inequality, item (iii.) from Lemma \ref{properties-mixed-morrey}.
\end{proof}
The next proposition is the corresponding Proposition 2.7 from \cite{Kozono-Yamazaki-semilinear-1994}.
\begin{prop}
Let us consider $\vec{q}\in[1,\infty)^{d}$, $\vec{\lambda}\in[0,1)^{d}$, $r\in [1,\infty]$ and $s\in\mathbb{R}$. Let $\left(u_{l}\right)_{l=-\infty}^{+\infty}$ be a sequence of tempered distributions on $\mathbb{R}^{d}$ such that $\supp\,{\cal F}u_{l}\subset {\cal A}_{l}$ and $u_{l}\in {\cal M}_{\vec{q},\vec{\lambda}}(\mathbb{R}^{d})$, for each $l\in \mathbb{Z}$. Let 
\begin{equation}
 B=\parallel 2^{sl}\parallel u_{l}\parallel_{\vec{q},\vec{\lambda}}\parallel_{l^{r}(\mathbb{Z})}\,\,\mbox{be finite}.
\end{equation}
Then, the sum $\sum_{l=0}^{\infty}u_{l}$ converge
\begin{itemize}
\item [(i.)] in $\dot{\cal N}^{s}_{\vec{q},\vec{\lambda},r}(\mathbb{R}^{d})$, if $r<\infty$, and
\item [(ii.)] in ${\cal N}^{s+\varepsilon}_{\vec{q},\vec{\lambda},1}(\mathbb{R}^{d})\cup {\cal N}^{s-\varepsilon}_{\vec{q},\vec{\lambda},1}(\mathbb{R}^{d})$, for each $\varepsilon>0$, if $r=\infty$. 
\end{itemize}
Moreover, if $u$ is the limit of the above convergent series, then
\begin{equation}
 \parallel u\parallel_{{\cal N}^{s}_{\vec{q},\vec{\lambda},r}(\mathbb{R}^{d})}\leq CB,
\end{equation}
where $C$ is a positive constant independent on $B$.
\end{prop}
\begin{proof}
The proof follows the same lines as in Theorem 2.7 from \cite{Kozono-Yamazaki-semilinear-1994}.   
\end{proof}

The next proposition is an extension of Theorem 2.9 from \cite{Kozono-Yamazaki-semilinear-1994}.
%%%%%%%%%%%%%%%%%%%%%%%%%%%%%%%%%%%%%%%
\begin{prop}
For $\gamma, s\in\mathbb{R}$, $\vec{q}\in[1,\infty)^{d}$, $\vec{\lambda}\in [0,1)^{d}$, there exist a positive integer $M$ and a positive constant $C$, such that the following estimate is valid for the operator $P(D)u={\cal F}^{-1}[P(\eta)\hat{u}]$, for each $B>0$.
\begin{itemize}
\item Let $P(\eta)$ be a $C^{M}$-function on $\mathbb{R}^{d}-\left\{0\right\}$ for which 
\begin{equation}
 \mid \left( \partial^{\mid \alpha\mid}/\partial\eta^{\alpha}\right)P(\eta)\mid\leq B\mid \eta\mid^{\gamma-\mid \alpha\mid}, 
 \end{equation}
 for each $\alpha\in\mathbb{N}^{d}$ such that $\mid \alpha\mid \leq M$. Then, the operator $P(D)$, defined on the space ${\cal S}'/{\cal P}$, is bounded from ${\cal N}^{s}_{\vec{q},\vec{\lambda},r}(\mathbb{R}^{d})$ to ${\cal N}^{s-\gamma}_{\vec{q},\vec{\lambda},r}(\mathbb{R}^{d})$   
and satisfies the estimate
\begin{equation}
\parallel P(D)u\parallel_{{\cal N}^{s-\gamma}_{\vec{q},\vec{\lambda},r}(\mathbb{R}^{d})}\leq CB
\parallel P(D)u\parallel_{{\cal N}^{s}_{\vec{q},\vec{\lambda},r}(\mathbb{R}^{d})}.
\end{equation}
\end{itemize}
\end{prop}
\begin{proof}
The proof follows the same lines as in Theorem 2.9 from \cite{Kozono-Yamazaki-semilinear-1994}.   
\end{proof}
The next proposition is an extension of Proposition 2.10 from \cite{Kozono-Yamazaki-semilinear-1994}.
\begin{prop}\label{realization}
Let $\vec{q}\in[1,\infty)^{d}$, $\vec{\lambda}\in[0,1)^{d}$, $r\in[1,\infty]$ and $s\in\mathbb{R}$. Suppose that 
\begin{equation}
 s<\displaystyle{\sum_{i=1}^{d}}\left(\frac{1-\lambda_{i}}{q_{i}}\right)\,\,\mbox{or}\,\, s=\displaystyle{\sum_{i=1}^{d}}\left(\frac{1-\lambda_{i}}{q_{i}}\right)\,\,\mbox{and}\,\,r=1.   
\end{equation}
Then, for each $u\in {\cal N}^{s}_{\vec{q},\vec{\lambda},r}(\mathbb{R}^{d})$, the sequence $(u_{(k)})_{k=0}^{\infty}$, given by
\begin{equation}
 u_{(k)}=\displaystyle{\sum_{l=-k}^{\infty}}\dot{\Delta}_{l}u,   
\end{equation}
converges in ${\cal S}'$, and the limit is the canonical representative of $u$ in ${\cal S}'/{\cal P}$.
\end{prop}
\begin{proof}
 The proof is a rewritten, in our context, from Proposition 2.10 from \cite{Kozono-Yamazaki-semilinear-1994}.   
\end{proof}
The next proposition is a partial extension of Proposition 2.11 from \cite{Kozono-Yamazaki-semilinear-1994}.
\begin{prop}\label{auxiliar-propo-inclusion}
Let $\vec{q}\in[1,\infty)^{d}$ and $\vec{\lambda}\in[0,1)^{d}$. Then, 
\begin{itemize}
\item  [(i.)] the continuous inclusions are valid
\begin{equation}
 {\cal N}^{0}_{\vec{q},\vec{\lambda}, 1}(\mathbb{R}^{d})\hookrightarrow {\cal M}_{\vec{q},\vec{\lambda}}(\mathbb{R}^{d})\hookrightarrow {\cal N}^{0}_{\vec{q},\vec{\lambda},\infty}(\mathbb{R}^{d}),\,\, \mbox{and}
\end{equation}
\item [(ii.)] the continuous inclusion is valid
\begin{equation}
 {\cal M}_{\vec{\lambda}}(\mathbb{R}^{d})\hookrightarrow {\cal N}^{0}_{\vec{1},\vec{\lambda}.
 \infty}(\mathbb{R}^{d}).    
\end{equation}
\end{itemize}
\end{prop}
\begin{proof}
 The proof follows the lines of Proposition 2.11 from \cite{Kozono-Yamazaki-semilinear-1994}.    
\end{proof}
Let us observe that in the notation of Besov-Morrey spaces from \cite{Kozono-Yamazaki-semilinear-1994}, for $\vec{\lambda}\in[0,1)^{d}$, we have ${\cal N}^{0}_{\vec{1},\vec{\lambda}.
 \infty}(\mathbb{R}^{d})={\cal N}^{0}_{p,1.
 \infty}(\mathbb{R}^{d})$, for $p=\frac{d}{\sum_{i=1}^{d}(1-\lambda_{i})}$.

The next Lemma is an extension of Lemma 2.4 from \cite{Bahouri-Chemin-Danchin-Fourier-PDE-2011}.
\begin{lema}
Let ${\cal A}$ be a given annulus. There exist two positive constants $c$ and $C$ such that, for each $\vec{q}\in[1,\infty)^{d}$, $\vec{\lambda}\in[0,1)^{d}$ and any pair of positive real numbers $t$ and $L$, we have
\begin{equation}
 \displaystyle{\sup\,\hat{u}}\subset L{\cal A}\,\,\mbox{implies}\,\,\parallel e^{\nu t\Delta}u\parallel_{{\cal M}_{\vec{q},\vec{\lambda}}(\mathbb{R}^{d})}\leq Ce^{-c\nu tL^{2}}\parallel u\parallel_{{\cal M}_{\vec{q},\vec{\lambda}}(\mathbb{R}^{d})},   
\end{equation}
for each $u\in{\cal M}_{\vec{q},\vec{\lambda}}(\mathbb{R}^{d})$ and each $t\in (0,\infty)$. Here, in both the itens, $\nu$ is a given positive real number.
\end{lema}
\begin{proof}
In order to prove this lemma we use the Young's inequality, item (iii.) from Lemma \ref{properties-mixed-morrey}, and follow the reasoning for the proof of Lemma 2.4 from \cite{Bahouri-Chemin-Danchin-Fourier-PDE-2011}.  
\end{proof}
%%%%%%%%%%%%%%%%%%%%%%%%%%%%%%%%%%%%%%%%%%%%%%%%
%%%%%%%%%%%%%%%%%%%%%%%%%%%%%%%%%%%%%%%%%%%%%

%%%%%%%%%%%%%%%%%%%%%%%%%%%%%%%%%%%%%%%%%%%%%%%%%%%
\section{Fourier-Besov mixed-Morrey spaces}\label{section-FBMMorrey}
Now, inspired on Fourier-Besov spaces and Fourier-Besov-Morrey spaces, we introduce the Fourier-Besov mixed-Morrey spaces, which are new in the study of Navier-Stokes equations.
\begin{defi}(Fourier-Besov mixed-Morrey spaces)

Let us consider four parameters 
$$
    \text{$\vec{q}\in[1,\infty)^{d}$, $\vec{\lambda}\in[0,1)^{d}$, $r\in[1,\infty]$ and $s\in\mathbb{R}$.}
$$ 
The $d$-dimensional homogeneous Fourier-Besov mixed-Morrey spaces are defined as
\begin{itemize}
\item [(i.)] For $\vec{q}\in[1,\infty)^{d}$, $\vec{\lambda
}\in[0,1)^{d}$, $r\in[1,\infty)$ and $s\in\mathbb{R}$,
\begin{equation}
\dot{{\cal FN}}^{s}_{\vec{q},\vec{\lambda},r}(\mathbb{R}^{d})=\lbrace
g\in {\cal S}'/{\cal P};\parallel g\parallel_{\dot{{\cal FN}}^{s}_{\vec{q},\vec{\lambda},r}(\mathbb{R}^{d})}=
\left(\displaystyle{\sum_{l\in\mathbb{Z}}2^{lsr}\parallel \phi_{l}\hat{g}\parallel^{r}_{L^{\vec{p}}(\mathbb{R}^{d})}}
\right)^{\frac{1}{r}}\,\,\mbox{is finite}
\rbrace.
\end{equation}
\item [(ii.)]  For  $\vec{q}\in[1,\infty)^{d}$, $\vec{\lambda}\in[0,1)^{d}$, $r=\infty$ and $s\in\mathbb{R}$,
\begin{equation}
\dot{{\cal FN}}^{s}_{\vec{q},\vec{\lambda},\infty}(\mathbb{R}^{d})=\lbrace
g\in {\cal S}'/{\cal P};\parallel g\parallel_{\dot{{\cal FN}}^{s}_{\vec{q},\vec{\lambda},\infty}(\mathbb{R}^{d})}=\displaystyle{\sup_{l\in\mathbb{Z}} 2^{ls}\parallel \phi_{l}\hat{g}\parallel_{{\cal M}_{\vec{q},\vec{\lambda}}(\mathbb{R}^{d})}}\,\,\mbox{is finite}
\rbrace. 
\end{equation}
\end{itemize}
We also called these spaces as the $d$-dimensional homogeneous mixed-norm Fourier-Besov-Morrey spaces or the $d$-dimensional homogeneous Fourier-Besov spaces based on Morrey spaces.
\end{defi}
%%%%%%%%%%%%%%%%%%%%%%%%%%%%%%%%%%%%%%%%%%%%%%%%
The next properties could be considered as an extension of properties in last lemma on Fourier-Besov-Morrey spaces.
\subsection{Proof of Theorem \ref{fourier-besov-properties}}
\begin{proof}
 \begin{itemize}
\item [(i.)] The proof os this item is based on the H\"older's inequality and the definition of mixed-Morrey spaces as follows,
\begin{equation}
\begin{split}
  \parallel \phi_{l}\hat{f}\parallel_{L^{1}(\mathbb{R}^{d})} &=
  \parallel \chi_{\{\xi;\frac{3}{4}\cdot 2^{l}\leq \mid\xi\mid\leq \frac{8}{3}\cdot 2^{l}\}}(\cdot) \phi_{l}\hat{f}\parallel_{L^{1}(\mathbb{R}^{d})}\\
  &\leq \parallel \chi_{\{\xi;\frac{3}{4}\cdot2^{l}\leq \mid\xi\mid\leq \frac{8}{3}\cdot2^{l}\}}(\cdot)\parallel_{L^{\vec{q}'}(\mathbb{R}^{d})}\parallel \phi_{l}\hat{f}\cdot\chi_{B(0,\frac{8}{3}\cdot2^{l})}\parallel_{L^{\vec{q}}(\mathbb{R}^{d})}\\
  &\leq C_{\vec{q}} 2^{l(\sum_{i=1}^{d}\frac{1}{q_{i}'})}\cdot 2^{l(\sum_{i=1}^{d})\frac{\lambda_{i}}{q_{i}}}\parallel \phi_{l}\hat{f}\parallel_{{\cal M}_{\vec{q},\vec{\lambda}}(\mathbb{R}^{d})}\\
  &\leq C_{\vec{q}}2^{l(\sum_{i=1}^{d}\frac{1}{q_{i}'} + \frac{\lambda_{i}}{q_{i}})}
  \parallel \phi_{l}\hat{f}\parallel_{{\cal M}_{\vec{q},\vec{\lambda}}(\mathbb{R}^{d})}
\end{split}
\end{equation}
which implies the desired estimate (i.).
\item[(ii)] In order to prove iten (ii.) it enough to consider the next H\"older's inequality in mixed-Morrey spaces, 
\begin{equation}
 \begin{split}
 \parallel \phi_{l}\hat{f}\parallel_{\vec{q},\vec{\lambda}}&\leq \parallel \chi_{\{\xi;\frac{3}{4}\cdot 2^{l}\leq \mid\xi\mid\leq \frac{8}{3}\cdot 2^{l} \}}(\cdot)\parallel_{\vec{p},\vec{\alpha}}\parallel \phi_{l}\hat{f}\parallel_{\vec{r},\vec{\mu}}\\
 &\leq C_{\vec{q},\vec{r}}2^{l(\sum_{i=1}^{d}\frac{1-\alpha_{i}}{p_{i}})}\parallel \phi_{l}\hat{f}\parallel_{\vec{r},\vec{\mu}}\\
 &=C_{\vec{q},\vec{r}}2^{l(\sum_{i=1}^{d}\frac{1}{q_{i}}-\frac{1}{r_{i}}-\frac{\lambda_{i}}{q_{i}}+\frac{\mu_{i}}{r_{i}})}\parallel \phi_{l}\hat{f}\parallel_{\vec{r},\vec{\mu}},
 \end{split}   
\end{equation}
which implies the desired estimate iten (ii.).
In fact, the assumptions $q_{i}\geq 1, r_{i}>q_{i}, \frac{\mu_{i}}{r_{i}}\leq \frac{\lambda_{i}}{q_{i}}$ and $\frac{1-\mu_{i}}{r_{i}}<\frac{1-\lambda_{i}}{q_{i}}$ implies that the choose 
\begin{equation}
 p_{i}=\frac{q_{i}r_{i}}{r_{i}-q_{i}}\,\, \mbox{and}\,\,\alpha_{i}=p_{i}\left(\frac{\lambda_{i}}{q_{i}}-\frac{\mu_{i}}{r_{i}}\right)
\end{equation}
satisfy $p_{1}\geq 1$, $p_{1}<\infty$, $\alpha_{i}\geq 0$ and $\alpha_{i}<1$, respectively.
\item [(iii.)] The proof of item (iii.1.) follows from item (ii) in Lemma \ref{properties-mixed-morrey} and to prove item (iii.2) it is enough to recall the continuous inclusion $l^{a}\hookrightarrow l^{b}$, for the given $a\leq b$.
\item[(iv.)] The proof of this item is based on the proved item (ii) above.
\item [(v.)] The proof follows from the proved item (i) and from the well-known estimate $$\parallel \dot{\Delta}_{l}u\parallel_{L^{\infty}(\mathbb{R}^{d})}\leq C\parallel \phi_{l}\hat{u}\parallel_{L^{1}(\mathbb{R}^{d})},$$ for some constant $C>0$.
\item [(vi.)] The proof of item (vi.) follows the lines from Theorem 2.4 from  \cite{Kozono-Yamazaki-semilinear-1994}.
\item [(vii.)] It follows from the proved itens (v.1) and (v.2) above from Theorem \ref{fourier-besov-properties}, since the well-known fact that $\dot{{\cal B}}^{\sigma}_{\infty,r}(\mathbb{R}^{d})\hookrightarrow {\cal S}'/{\cal P}$, for all $\sigma\in\mathbb{R}$ and $r\in[1,\infty]$. 
 
 \end{itemize}   
\end{proof}
  
%%%%%%%%%%%%%%%%%%%%%%%%%%%%%%%%%%%%%%%%%%%%%%%%%%
%%%%%%%%%%%%%%%%%%%%%%%%%%%%%%%%%%%%%%%%%%%%%%%%%

%%%%%%%%%%%%%%%%%

%%%%%%%%%%%%%%%%%%%%%%%%%%%%%%
\section{Application to Navier-Stokes system}\label{Section-NS-system}
%The main tools for the method and linear estimates}
\label{main-tools}
%%%%%%%%%%%%%%%%%%%%%%%%%%%%%%
In this section we use the Bony's paraproduct approach and the introduced Besov mixed-Morrey spaces and Fourier-Besov mixed-Morrey spaces in order to get a global well-posedness result for the Navier-Stokes system. In the first subsection, we first recall the Stokes semigroup and the definition of mild solutions. In the second subsection we set the main linear and bilinear estimates. In the last subsection we sketch the proof of the last two main theorems.

\subsection{The Stokes semigroup and Mild solutions}
We shall denote by $S_{\nu}(t)=e^{\nu t\Delta}$, for $t>0$, the Stokes semigroup whose Fourier representation is given by 
\begin{equation*}
    [S_{\nu}(t)g]^{\wedge}(\eta)=e^{-\nu t\mid \eta\mid^{2}}\hat{g}(\eta).  
\end{equation*}
Let us denote by $\mathbb{P}$ the Leray projector, which in Fourier variables is given by its symbols $\hat{\mathbb{P}}$, that is, 
\begin{equation*}
 \left(\hat{\mathbb{P}}(\eta)\right)_{i,j}=\delta_{i,j}-\frac{\eta_{i}\eta_{j}}{\mid \eta\mid^{2}}, \,\,\mbox{for each}\,\, i,j=1,\dots, d,   
\end{equation*}
where $\delta_{i,j}$ denotes the Kronecker delta.

The mild formulation for the Navier-Stokes equations is given by the next abstract equation
\begin{equation}\label{mild-formulation}
u(t)=S_{\nu}(t)u_{0}+B(u,u)(t),    
\end{equation}
for all $t\in[0,\infty)$, where the bilinear operator $B(\cdot,\cdot)$ is given by 
\begin{equation*}
 B(v,w)(t)=\displaystyle{\int_{0}^{t}}S_{\nu}(t-t')\mathbb{P}\mbox{div}\left(v\otimes w\right)(t')\,dt',   
\end{equation*}
which is given in Fourier variables by 
\begin{equation*}
[B(v,w)]^{\wedge}(\eta,t)=\displaystyle{\int_{0}^{t}}e^{-\nu(t-t')\mid\eta\mid^{2}}\hat{\mathbb{P}}(\eta)[i\eta\cdot \left(v\otimes w\right)^{\wedge}(\eta,t')]\,dt', 
\end{equation*}
%{\bf acima nao faz sentido produto escalar, $\eta$-vetor e $(u\otimes v)^{\wedge}$-matriz} 
where the components of the matrix $(u\otimes v)^{\wedge}$ and the vector $i\eta\cdot\left(v\otimes w\right)^{\wedge}$ are given by 
\begin{equation*}
  \left(v\otimes w\right)^{\wedge}_{i,j}=\hat{v}_{i}\ast\hat{w}_{j}, \,\,\mbox{for}\,\,i,j=1,\dots,d,\,\,\mbox{and}\,\,\left(i\eta\cdot\left(v\otimes w\right)^{\wedge}\right)_{l}=i\eta\cdot\mbox{row}_{l}\,\left(v\otimes w\right)^{\wedge},\,\,\mbox{for}\,\,l=1,\dots,d.  
\end{equation*}
\begin{defi}
 Let $Z$ be a suitable Banach space. We say that $u\in Z$ is a mild solution for the Navier-Stokes equations if $u$ satisfies the mild formulation (\ref{mild-formulation}) above.
 \begin{obs}
In the context of the main theorems below $Z={\cal L}^{\infty}\left(I;{\dot{{\cal N}}^{\sigma}_{\vec{q},\vec{\lambda},r}(\mathbb{R}^{d})}\right)\bigcap {\cal L}^{1}\left(I;{\dot{{\cal N}}^{\sigma+2}_{\vec{q},\vec{\lambda},r}(\mathbb{R}^{d})}\right)$, so that in order to get $S_{\nu}(\cdot)u_{0}\in Z$, we required the initial data $u_{0}\in \dot{{\cal N}}^{\sigma}_{\vec{q},\vec{\lambda},r}(\mathbb{R}^{d})$. Similar remark for Fourier-Besov mixed-Morrey spaces.  
 \end{obs}
 %{\bf com o dado inicial dado no mesmo espaco, importante ressaltar na definicao mild }
\end{defi}
Let us recall the space ${\cal L}^{a}(I;\dot{{\cal FN}}^{s}_{\vec{q},\vec{\lambda},r}(\mathbb{R}^{d}))$ of Bochner measurable functions defined on interval $I$ toward $\dot{{\cal FN}}^{s}_{\vec{q},\vec{\lambda},r}(\mathbb{R}^{d})$ provided by the norm
\begin{equation*}
\parallel g\parallel_{{\cal L}^{a}(I;\dot{{\cal FN}}^{s}_{\vec{q},\vec{\lambda},r}(\mathbb{R}^{d}))}=\displaystyle{\parallel
2^{js} \parallel\phi_{l}\hat{g}\parallel_{L^{a}\left(I;{\cal M}_{\vec{q},\vec{\lambda}}(\mathbb{R}^{d})\right)}\parallel_{l^{r}(\mathbb{Z})}}.    
\end{equation*}
A similar definition should be considered if we consider $\dot{{\cal N}}^{\sigma}_{\vec{q},\vec{\lambda},r}(\mathbb{R}^{d})$ instead of $\dot{{\cal FN}}^{s}_{\vec{q},\vec{\lambda},r}(\mathbb{R}^{d})$.

%%%%%%%%%%%%%%%%%%%
\subsection{Linear estimates and Bilinear Estimates}

In this subsection we established linear estimates involving the Stokes semigroup, 
which we shall use for computations in relation to the initial data and also for the auxiliary linear operator related to the bilinear term.

\subsubsection{The linear estimates}
Here, we proved the main linear estimates related to Besov mixed-Morrey spaces and Fourier-Besov mixed-Morrey spaces, involving the Stokes semigroup.
\begin{lema}\label{estimates-lin-besov}
 Let us take $\vec{q}\in[1,\infty)^{d}$, $\vec{\lambda}\in[0,1)^{d}$, $r\in[1,\infty]$ and $\sigma\in\mathbb{R}$. There exists a constant $C>0$, independent on $\nu$, such that
\begin{equation*}
\begin{split}
 (i.) &\parallel S_{\nu}(\cdot)u_{0}\parallel_{{\cal L}^{\infty}\left(I;\dot{{\cal N}}^{\sigma}_{\vec{q},\vec{\lambda},r}(\mathbb{R}^{d})\right)}\leq C \parallel u_{0}\parallel_{\dot{{\cal N}}^{\sigma}_{\vec{q},\vec{\lambda},r}(\mathbb{R}^{d})}, 
 \\
 (ii.)&\parallel S_{\nu}(\cdot)u_{0}\parallel_{{\cal L}^{1}\left(I;\dot{{\cal N}}^{\sigma+2}_{\vec{q},\vec{\lambda},r}(\mathbb{R}^{d})\right)}\leq \frac{C}{\nu} \parallel u_{0}\parallel_{\dot{{\cal N}}^{\sigma}_{\vec{q},\vec{\lambda},r}(\mathbb{R}^{d})},  
\end{split}
\end{equation*}
 for any $u_{0}\in \dot{{\cal N}}^{\sigma}_{\vec{q},\vec{\lambda},r}(\mathbb{R}^{d})$.
\end{lema}
\begin{proof}
\begin{itemize}
\item [1.] Proof of item (i.): From $e^{-ct2^{2l}}\leq 1$, for all $t>0$ and any $l\in\mathbb{Z}$ we obtain (i.) if we consider the definition of the norm on ${\cal L}^{\infty}\left(I;\dot{{\cal N}}^{\sigma}_{\vec{q},\vec{\lambda},r}(\mathbb{R}^{d})\right)$. 
\item [2.] Proof of item (ii.): Since $\parallel e^{\nu t\Delta}\dot{\Delta}_{l}u\parallel_{{\cal M}_{\vec{q},\vec{\lambda}}(\mathbb{R}^{d})}\leq C e^{-c\nu t2^{2l}} \parallel\dot{\Delta}_{l}u\parallel_{{\cal M}_{\vec{q},\vec{\lambda}(\mathbb{R}^{d})}}$, for some constants $c$ and $C$, independent on $l$ and $u$, and the definition of the norm on  ${\cal L}^{1}\left(I;\dot{{\cal N}}^{\sigma+2}_{\vec{q},\vec{\lambda},r}(\mathbb{R}^{d})\right)$, we get
\begin{equation*}
 \begin{split}
 \parallel S_{\nu}u_{0}\parallel_{{\cal L}^{1}\left(I;\dot{{\cal N}}^{\sigma+2}_{\vec{q},\vec{\lambda},r}(\mathbb{R}^{d})\right)}&\leq
 \parallel 2^{l(\sigma+2)}\displaystyle{\int_{0}^{\infty}}\parallel e^{\nu t'\Delta}\dot{\Delta}_{l}u_{0}\parallel_{{\cal M}_{\vec{q},\vec{\lambda}}(\mathbb{R}^{d})}\,dt'\parallel_{l^{r}(\mathbb{Z})}\\
 &\leq
 \parallel 2^{l(\sigma+2)}\left(\displaystyle{\int_{0}^{\infty}}e^{-c\nu t2^{2l}}\,dt'\right)\parallel\dot{\Delta}_{l}u_{0}\parallel_{{\cal M}_{\vec{q},\vec{\lambda}}(\mathbb{R}^{d})}\parallel_{l^{r}(\mathbb{Z})},
 \end{split}   
\end{equation*}
which implies (ii.) since the integral is bounded by $\frac{C}{\nu}2^{-2l}$.
\end{itemize}
\end{proof}

%%%%%%%%%%%%%%%%%%%%%%%%%%
\begin{lema}
Let us consider $\vec{q}\in[1,\infty)^{d}$, $\vec{\lambda}\in[0,1)^{d}$, $r\in[1,\infty]$ and $s\in \mathbb{R}$. There is a constant $C>0$, independent of $\nu$, such that
\begin{equation*}
\begin{split}
(i.)&\parallel S_{\nu}(\cdot)u_{0}\parallel_{{\cal L}^{\infty}\left(I;\dot{{\cal FN}}^{s}_{\vec{q},\vec{\lambda},r}(\mathbb{R}^{d})\right)}\leq C\parallel u_{0}\parallel_{\dot{{\cal FN}}^{s}_{\vec{q},\vec{\lambda},r}(\mathbb{R}^{d})},\\    
(ii.)&
\parallel S_{\nu}(\cdot)u_{0}\parallel_{{\cal L}^{1}\left(I;\dot{{\cal FN}}^{s+2}_{\vec{q},\vec{\lambda},r}(\mathbb{R}^{d})\right)}\leq \frac{C}{\nu}\parallel u_{0}\parallel_{\dot{{\cal FN}}^{s}_{\vec{q},\vec{\lambda},r}(\mathbb{R}^{d})},    
\end{split}
\end{equation*} 
for every $u_{0}\in \dot{{\cal FN}}^{s}_{\vec{q},\vec{\lambda},r}(\mathbb{R}^{d})$.
\end{lema}

\begin{proof}
\begin{itemize}
\item [1.]Proof of item (i.): Since $e^{-\nu t\mid\eta\mid^{2}}\leq 1$, for all $t\geq 0$ and $\eta\in\mathbb{R}^{d}$, we easily get $(i)$, by considering the definition of ${\cal L}^{\infty}\left(I;\dot{{\cal FN}}^{s}_{\vec{q},\vec{\lambda},r}(\mathbb{R}^{d})\right)$.  
\item [2.] Proof of item (ii.): By definition of the norm on ${\cal L}^{1}\left(I;\dot{{\cal FN}}^{s+2}_{\vec{q},\vec{\lambda},r}(\mathbb{R}^{d})\right)$, we get
\begin{equation*}
 \begin{split}
  \parallel S_{\nu}(\cdot)u_{0}\parallel_{{\cal L}^{1}\left(I;\dot{{\cal FN}}^{s+2}_{\vec{q},\vec{\lambda},r}(\mathbb{R}^{d})\right)}&\leq 
  \displaystyle{\parallel 
  2^{l(s+2)}\int_{0}^{\infty}\parallel e^{-\nu t'\mid \eta\mid^{2}}\phi_{l}(\eta)\hat{u}_{0}\parallel_{{\cal M}_{\vec{q},\vec{\lambda}}(\mathbb{R}^{d})}\,dt'
  \parallel}_{l^{r}(\mathbb{Z})}\\
  &\leq 
  \displaystyle{\parallel 
  2^{l(s+2)}\left(\int_{0}^{\infty} e^{-\nu t'2^{2l}}\,dt'\right)\parallel \phi_{l}(\eta)\hat{u}_{0}\parallel_{{\cal M}_{\vec{q},\vec{\lambda}}(\mathbb{R}^{d})}
  \parallel}_{l^{r}(\mathbb{Z})}.
 \end{split}   
\end{equation*}
Since the last integral is bounded by $\frac{C
}{\nu}2^{-2l}$, we have the desired $(ii.)$.
\end{itemize}
\end{proof}

%%%%%%%%%%%%%%%%%%%%%%%%%%%
\subsubsection{The bilinear estimates}
\label{bilinear-estimates-s}
%%%%%%%%%%%%%%%%%%%%%%%%%%%

In this section, we first introduce an auxiliary linear operator and we establish 
some basic estimates for this operator. We also set the precise estimate for 
the bilinear term associated to the mild formulation. 

\medskip
First, let us define the following auxiliary linear operator
\begin{equation}\label{auxiliar-linearop}
 A(g)(t,\cdot):= \displaystyle{\int_{0}^{t}}e^{\nu(t-t')\Delta}\mathbb{P}g(t',\cdot)\,dt'   
\end{equation}
which in Fourier variables is given by
\begin{equation*}
 [A(g)]^{\wedge}(\eta,t)=\displaystyle{\int_{0}^{t}e^{-\nu(t-t')\mid\eta\mid^{2}}\hat{\mathbb{P}}(\eta)\hat{g}(\eta,t')\,dt'}.   
\end{equation*}
Then, we have the following 
\begin{lema}
Let us take $\vec{q}\in[1,\infty)^{d}$, $\vec{\lambda}\in[0,1)^{d}$, $r\in[1,\infty]$ and $\sigma\in\mathbb{R}$. Then, there exists $C'>0$ independent on $\nu$ such that
\begin{equation*}
\begin{split}
    (i.)&\parallel A(g)\parallel_{{\cal L}^{\infty}\left(I;\dot{{\cal N}}^{\sigma}_{\vec{q},\vec{\lambda},r}(\mathbb{R}^{d})\right)}\leq C'
    \parallel g\parallel_{{\cal L}^{1}\left(I;\dot{{\cal N}}^{\sigma}_{\vec{q},\vec{\lambda},r}(\mathbb{R}^{d})\right)},
    \\
(ii.)& \parallel A(g)\parallel_{{\cal L}^{1}\left(I;\dot{{\cal N}}^{\sigma+2}_{\vec{q},\vec{\lambda},r}(\mathbb{R}^{d})\right)}\leq \frac{C'}{\nu}
    \parallel g\parallel_{{\cal L}^{1}\left(I;\dot{{\cal N}}^{\sigma}_{\vec{q},\vec{\lambda},r}(\mathbb{R}^{d})\right)},
\end{split}
\end{equation*}
for all $g\in {\cal L}^{1}\left(I;\dot{{\cal N}}^{\sigma}_{\vec{q},\vec{\lambda},r}(\mathbb{R}^{d})\right)$.
\end{lema}
\begin{proof}

First, since the Riesz transform ${\cal R}_{j}$ is given by ${\cal R}_{j}[\dot{\Delta}_{l}g]={\cal F}^{-1}[\frac{i\eta_{j}}{\mid \eta\mid} (\dot{\Delta}_{l}g)^{\wedge}(\eta)]$, due to the Young's inequality in ${\cal M}_{\vec{q},\vec{\lambda}}(\mathbb{R}^{d})$, we conclude that 
    \begin{equation*}
     \parallel {\cal R}_{j}\dot{{\Delta}}_{l}g\parallel_{{\cal M}_{\vec{q},\vec{\lambda}}(\mathbb{R}^{d})}\leq C \parallel \dot{{\Delta}}_{l}g\parallel_{{\cal M}_{\vec{q},\vec{\lambda}}(\mathbb{R}^{d})},  
    \end{equation*}
    for some constant $C>0$ independent on $j, l$ and $g$. Therefore, by definition of the Leray operator $\mathbb{P}$, we obtain
    \begin{equation*}
     \parallel \mathbb{P}\dot{{\Delta}}_{l}g\parallel_{{\cal M}_{\vec{q},\vec{\lambda}}(\mathbb{R}^{d})}\leq C \parallel\dot{{\Delta}}_{l}g\parallel_{{\cal M}_{\vec{q},\vec{\lambda}}(\mathbb{R}^{d})},  
    \end{equation*}
    for some constant $C>0$ independent on $l$ and $g$.
\begin{itemize}
    \item [1.] Proof of item (i.): Since $\parallel e^{\nu t\Delta}\dot{\Delta}_{l}u\parallel_{{\cal M}_{\vec{q},\vec{\lambda}}(\mathbb{R}^{d})}\leq e^{-c\nu t2^{2l}}\parallel \dot{\Delta}_{l}u\parallel_{{\cal M}_{\vec{q},\vec{\lambda}}(\mathbb{R}^{d})}$ and from the definition of the norm on  ${\cal L}^{\infty}\left(I;\dot{{\cal N}}^{\sigma}_{\vec{q},\vec{\lambda},r}(\mathbb{R}^{d})\right)$, we get
\begin{equation*}
\begin{split}
   \parallel A(g)\parallel_{{\cal L}^{\infty}\left(I;\dot{{\cal N}}^{\sigma}_{\vec{q},\vec{\lambda},r}(\mathbb{R}^{d})\right)}&\leq \parallel 2^{l\sigma}\parallel \displaystyle{\int_{0}^{t}}e^{-c\nu (t-t')2^{2l}}\parallel \dot{\Delta}_{l}g\parallel_{{\cal M}_{\vec{q},\vec{\lambda}}(\mathbb{R}^{d})}\,dt'\parallel_{L^{\infty}(I)}\parallel_{l^{r}(\mathbb{Z})}\\
   &\leq \parallel 2^{l\sigma}\parallel \parallel\dot{\Delta}_{l}g\parallel_{L^{1}\left((0,t);{\cal M}_{\vec{q},\vec{\lambda}}(\mathbb{R}^{d})\right)}\parallel_{L^{\infty}(I)}\parallel_{l^{r}(\mathbb{Z})}
\end{split}    
\end{equation*}
and this implies (i.) since $(0,t)\subset I$.
\item [2.] Proof of item (ii.): Again, by the definition of the norm on ${\cal L}^{1}\left(I;\dot{{\cal N}}^{\sigma+2}_{\vec{q},\vec{\lambda},r}(\mathbb{R}^{d})\right)$, we have 
\begin{equation*}
 \begin{split}
  \parallel A(g)\parallel_{{\cal L}^{1}\left(I;\dot{{\cal N}}^{\sigma+2}_{\vec{q},\vec{\lambda},r}(\mathbb{R}^{d})\right)}&\leq 
  \parallel 2^{l(\sigma+2)}\parallel \displaystyle{\int_{0}^{t}}e^{-c\nu(t-t')2^{2l}}\parallel \dot{\Delta}_{l}g\parallel_{{\cal M}_{\vec{q},\vec{\lambda}}(\mathbb{R}^{d})}\,dt'\parallel_{L^{1}(I)}\parallel_{l^{r}(\mathbb{Z})}\\
  &\leq \parallel 2^{l(\sigma+2)}\displaystyle{\int_{0}^{\infty}}
  \left(\displaystyle{\int_{t'}^{\infty}}e^{-c\nu(t-t')2^{2l}}\,dt\right)\cdot \parallel \dot{\Delta}_{l}g(t')\parallel_{{\cal M}_{\vec{q},\vec{\lambda}}(\mathbb{R}^{d})}\,dt'\parallel_{l^{r}(\mathbb{Z})}
 \end{split}   
\end{equation*}
and this implies (ii.) since the integral is bounded by $\frac{C}{\nu}2^{-2l}$.
\end{itemize}
\end{proof}
\begin{obs}
One observes that, 
the Riesz transform is bounded on homogeneous Besov mixed-Morrey spaces. 
In fact, the proof is similar to the usual homogeneous Besov spaces (\cite{Sawano-homogeneous-2020}, Theorem 1.3). 
This assertion implies, by definition, that the Leray projector $\mathbb{P}$ is also bounded on these homogeneous Besov mixed-Morrey spaces.   
\end{obs}
\begin{lema}
Let us consider $\vec{q}\in[1,\infty)^{d}$, $\vec{\lambda}\in[0,1)^{d}$, $r\in[1,\infty]$ and $s\in\mathbb{R}$. Then, there exists $C'>0$ independent on $\nu$ such that
\begin{equation*}
\begin{split}
(i.)&\parallel A(g)\parallel_{{\cal L}^{\infty}\left(I;\dot{{\cal FN}}^{s}_{\vec{q},\vec{\lambda},r}(\mathbb{R}^{d})\right)}\leq C'  
\parallel g\parallel_{{\cal L}^{1}\left((I;\dot{{\cal FN}}^{s}_{\vec{q},\vec{\lambda},r}(\mathbb{R}^{d})\right)},\\
(ii.)&\parallel A(g)\parallel_{{\cal L}^{1}\left(I;\dot{{\cal FN}}^{s+2}_{\vec{q},\vec{\lambda},r}(\mathbb{R}^{d})\right)}\leq \frac{C'}{\nu}  \cdot
\parallel g\parallel_{{\cal L}^{1}\left(I;\dot{{\cal FN}}^{s}_{\vec{q},\vec{\lambda},r}(\mathbb{R}^{d})\right)},    
\end{split}
\end{equation*}
for all $g\in {\cal L}^{1}\left(I;\dot{{\cal FN}}^{s}_{\vec{q},\vec{\lambda},r}(\mathbb{R}^{d})\right)$.
\end{lema}
\begin{proof}
First, recall that $\supp \phi_{l}\subset 2^{l}{\cal A}=\{\eta\in\mathbb{R}^{d};\frac{3}{4}\cdot 2^{l}\leq \mid \eta\mid\leq \frac{8}{3}\cdot 2^{l}\}$ and ${\cal R}_{j}[\dot{\Delta}_{l}g]={\cal F}^{-1}[\frac{i\eta_{j}}{\mid \eta\mid}\phi_{l}\hat{g}]$. Then $\mid [{\cal R}_{j}(\dot{\Delta}_{l}g)]^{\wedge}(\eta)\mid=\mid \frac{i\eta_{j}}{\mid\eta\mid}\mid \mid\phi_{l}\hat{g}\mid\leq \mid\phi_{l}\hat{g}\mid$, for all $\eta\in\mathbb{R}^{d}-\{0\}$. This implies that $\mid {\cal F}[\mathbb{P}\dot{\Delta}_{l}g](\eta)\mid\leq \mid \phi_{l}\hat{g}(\eta)\mid$, for all $\eta\in\mathbb{R}^{d}-\{0\}$.
%First we observe that $\parallel \hat{\mathbb{P}}(\eta)\parallel\leq 2$, for all $\eta\in\mathbb{R}^{d}-\{0\}$ and $\supp \phi_{l}\subset 2^{l}{\cal A}=\{\eta\in\mathbb{R}^{d};\frac{3}{4}\cdot 2^{l}\leq \mid \eta\mid\leq \frac{8}{3}\cdot 2^{l}\}$.
%{\bf aqui esta mal, melhorar esta escrita, digo acima.}
\begin{itemize}
\item [1.] Proof of item (i.): We note that, by definition of the norm on ${\cal L}^{\infty}\left(I;\dot{{\cal FN}}^{s}_{\vec{q},\vec{\lambda},r}(\mathbb{R}^{d})\right)$, we have
\begin{equation*}
\begin{split}
\parallel A(g)\parallel_{{\cal L}^{\infty}\left(I;\dot{{\cal FN}}^{s}_{\vec{q},\vec{\lambda},r}(\mathbb{R}^{d})\right)}&\leq
\parallel 2^{ls}\parallel \displaystyle{\int_{0}^{t}}
e^{-\nu(t-t')2^{2l}}\parallel \phi_{l}\hat{g}\parallel_{{\cal M}_{\vec{q},\vec{\lambda}}(\mathbb{R}^{d})}\,dt'\parallel_{L^{\infty}(I)}\parallel_{l^{r}(\mathbb{Z})}\\
&\leq \parallel 2^{ls}\parallel \parallel \phi_{l}\hat{g}\parallel_{L^{1}\left((0,t);{\cal M}_{\vec{q},\vec{\lambda}}(\mathbb{R}^{d})\right)}\parallel_{L^{\infty}(I)}\parallel_{l^{r}(\mathbb{Z})}
\end{split}    
\end{equation*}
which implies iten (i.) since $(0,t)\subset I$.
\item [2.] Proof of item (ii.): Again, by definition of the norm on ${\cal L}^{1}\left(I;\dot{{\cal FN}}^{s+2}_{\vec{q},\vec{\lambda},r}(\mathbb{R}^{d})\right)$, we have
\begin{equation*}
\begin{split}
\parallel A(g)\parallel_{{\cal L}^{1}\left(I;\dot{{\cal FN}}^{s+2}_{\vec{q},\vec{\lambda},r}(\mathbb{R}^{d})\right)}&\leq
\parallel 2^{l(s+2)}\parallel \displaystyle{\int_{0}^{t}}
e^{-\nu(t-t')2^{2l}}\parallel \phi_{l}\hat{g}\parallel_{{\cal M}_{\vec{q},\vec{\lambda}}(\mathbb{R}^{d})}\,dt'\parallel_{L^{1}(I)}\parallel_{l^{r}(\mathbb{Z})}\\
&\leq \parallel 2^{l(s+2)}\cdot\frac{C'}{\nu}\cdot 2^{-2l} \parallel \phi_{l}\hat{g}\parallel_{L^{1}\left(I;{\cal M}_{\vec{q},\vec{\lambda}}(\mathbb{R}^{d})\right)}\parallel_{l^{r}(\mathbb{Z})}
\end{split}    
\end{equation*}
which implies iten (ii.), by considering the sum of the exponents of $2$.
\end{itemize}
\end{proof}
%%%%%%%%%%%%%%%%%%%%%%%%%%%%%%%%%%%%
\begin{lema}\label{estimates-bilinear-besov}
 Let us consider $\vec{q}\in[1,\infty)^{d}$ and $\vec{\lambda}\in[0,1)^{d}$ be such $\displaystyle{\sum_{i=1}^{d}}\frac{1-\lambda_{i}}{q_{i}}>0$, $r\in[1,\infty]$ and $\sigma=-1+\displaystyle{\sum_{i=1}^{d}}\frac{1-\lambda_{i}}{q_{i}}$. Then, there exists a constant $K_{0}>0$, independent on $\nu$, such that
 \begin{equation*}
  \parallel B(v,w)\parallel_{Z}\leq K_{0}\max{(1,\frac{1}{\nu})}\parallel v\parallel_{Z}\parallel w\parallel_{Z},   
 \end{equation*}
 for each $v,w\in Z={\cal L}^{\infty}\left(I;\dot{{\cal N}}^{\sigma}_{\vec{q},\vec{\lambda},r}(\mathbb{R}^{d})\right)\bigcap {\cal L}^{1}\left(I;\dot{{\cal N}}^{\sigma+2}_{\vec{q},\vec{\lambda},r}(\mathbb{R}^{d})\right)$.
\end{lema}
\begin{proof}
First, we observe that since $B(v,w)(t,\cdot)=A(\mbox{div}(v\otimes w))(t,\cdot)$, where the linear auxiliar operator $A(\cdot)$ is given in (\ref{auxiliar-linearop}), we get 
%{\bf aqui nao sabemos quem eh $A$, esta definido la traz: operador linear auxiliar...}
$$
\begin{aligned}
 \parallel B(v,w)\parallel_{Z}&\leq c\max{(1,\frac{1}{\nu})} \parallel \mbox{div}\left(v\otimes w\right)\parallel_{{\cal L}^{1}\left(I;\dot{{\cal N}}^{\sigma}_{\vec{q},\vec{\lambda},r}(\mathbb{R}^{d})\right)}
 \\
& \leq c\max{(1,\frac{1}{\nu})}\parallel v\otimes w\parallel_{{\cal L}^{1}\left(I;\dot{{\cal N}}^{\sigma+1}_{\vec{q},\vec{\lambda},r}(\mathbb{R}^{d})\right)},   
\end{aligned}
$$
where $c$ is a constant which does not depend on $\nu$. Therefore, it is enough to prove that 
\begin{equation*}
 \parallel v\otimes w\parallel_{{\cal L}^{1}\left(I;\dot{{\cal N}}^{\sigma+1}_{\vec{q},\vec{\lambda},r}(\mathbb{R}^{d})\right)}\leq \parallel v\parallel_{Z}\parallel w\parallel_{Z}.
\end{equation*}
In order to prove this estimate we consider the Bony's paraproduct decomposition
\begin{equation*}
 \begin{split}
  \dot{\Delta}_{l}\left(vw\right)&=\displaystyle{\sum_{\mid l-l'\mid\leq 4}}\dot{\Delta}_{l}\left( \dot{S}_{l'-1}v\dot{\Delta}_{l'}w\right) + \displaystyle{\sum_{\mid l-l'\mid\leq 4}}\dot{\Delta}_{l}\left( \dot{S}_{l'-1}w\dot{\Delta}_{l'}v\right) + \displaystyle{\sum_{l'\geq l-2}}\dot{\Delta}_{l}\left(\dot{\Delta}_{l'}v\tilde{\dot{\Delta}}_{l'}w\right)\\
  &=I^{l}_{1} + I^{l}_{2} + I^{l}_{3}.
 \end{split}   
\end{equation*}
Estimates for $I^{l}_{1}$: We first observe that the Bernstein's type inequality Lemma 2.4 (i.) (\ref{properties-Besov-spaces-mixed}) implies $\parallel \dot{\Delta}_{l''}v\parallel_{L^{\infty}(\mathbb{R}^{d})}\leq 2^{l''(\sigma +1)}\parallel \dot{\Delta}_{l''}v\parallel_{{\cal M}_{\vec{q},\vec{\lambda}}(\mathbb{R}^{d})}$ and thus we have
\begin{equation*}
\begin{split}
\parallel I^{l}_{1}\parallel_{L^{1}\left(I;{\cal M}_{\vec{q},\vec{\lambda}}(\mathbb{R}^{d})\right)}&\leq 
\displaystyle{\sum_{\mid l-l'\mid\leq 4}}\parallel \parallel \dot{S}_{l'-1}v\dot{\Delta}_{l'}w\parallel_{{\cal M}_{\vec{q},\vec{\lambda}}(\mathbb{R}^{d})}\parallel_{L^{1}(I)}\\
&\leq \displaystyle{\sum_{\mid l-l'\mid\leq 4}}
\parallel \displaystyle{\sum_{l''\leq l'}}\parallel\dot{\Delta}_{l''}v\parallel_{L^{\infty}(\mathbb{R}^{d})}\parallel\dot{\Delta}_{l'}w\parallel_{{\cal M}_{\vec{q},\vec{\lambda}}(\mathbb{R}^{d})}\parallel_{L^{1}(I)}\\
&\leq \displaystyle{\sum_{\mid l-l'\mid\leq 4}} \left(\displaystyle{\sum_{l''\leq l'}}
2^{l''(\sigma +1)}\parallel \dot{\Delta}_{l''}v\parallel_{L^{\infty}\left(I;{\cal M}_{\vec{q},\vec{\lambda}}(\mathbb{R}^{d})\right)}\right) \parallel \dot{\Delta}_{l'}w\parallel_{L^{1}\left(I;{\cal M}_{\vec{q},\vec{\lambda}}(\mathbb{R}^{d})\right)}\\
&= \displaystyle{\sum_{\mid l-l'\mid\leq 4}} \left(\displaystyle{\sum_{l''\leq l'}}
2^{l''}\cdot 2^{l''\sigma}\parallel \dot{\Delta}_{l''}v\parallel_{L^{\infty}\left(I;{\cal M}_{\vec{q},\vec{\lambda}}(\mathbb{R}^{d})\right)}\right) \parallel \dot{\Delta}_{l'}w\parallel_{L^{1}\left(I;{\cal M}_{\vec{q},\vec{\lambda}}(\mathbb{R}^{d})\right)}.
\end{split}   
\end{equation*}
By H\"older's inequality for series, we get
\begin{equation*}
 \begin{split}
 \parallel I^{l}_{1}\parallel_{L^{1}\left(I;{\cal M}_{\vec{q},\vec{\lambda}}(\mathbb{R}^{d})\right)}&\leq    
  \displaystyle{\sum_{\mid l-l'\mid\leq 4}}\left(\displaystyle{\sum_{l''\leq l'}2^{l''r'}}\right)^{\frac{1}{r'}}\cdot \parallel \dot{\Delta}_{l'}w\parallel_{L^{1}\left(I;{\cal M}_{\vec{q},\vec{\lambda}}(\mathbb{R}^{d})\right)}\parallel v \parallel_{{\cal L}^{\infty}\left(I;\dot{{\cal N}}^{\sigma}_{\vec{q},\vec{\lambda},r}(\mathbb{R}^{d})\right)}\\
  &\leq \displaystyle{\sum_{\mid l-l'\mid\leq 4}}2^{l'} \parallel \dot{\Delta}_{l'}w\parallel_{L^{1}\left(I;{\cal M}_{\vec{q},\vec{\lambda}}(\mathbb{R}^{d})\right)}\parallel v\parallel_{Z}\\
&=  \displaystyle{\sum_{\mid l-l'\mid\leq 4}}2^{-l'(\sigma+1)}\cdot 2^{l'(\sigma+2)}\parallel \dot{\Delta}_{l'}w\parallel_{L^{1}\left(I;{\cal M}_{\vec{q},\vec{\lambda}}(\mathbb{R}^{d})\right)}\parallel v\parallel_{Z}.
 \end{split}   
\end{equation*}
Multiplying by $2^{l(\sigma+1)}$ in both sides, we have
\begin{equation*}
 \begin{split}
  2^{l(\sigma+1)}\parallel I^{l}_{1}\parallel_{L^{1}\left(I;{\cal M}_{\vec{q},\vec{\lambda}}(\mathbb{R}^{d})\right)}&\leq
  \displaystyle{\sum_{l'}}2^{(l-l')(\sigma+1)}\chi_{\{k;\mid k\mid\leq 4\}}(l-l')\cdot \left(2^{l'(\sigma+2)}\parallel \dot{\Delta}_{l'}w\parallel_{L^{1}\left(I;{\cal M}_{\vec{q},\vec{\lambda}}(\mathbb{R}^{d})\right)}\right)\cdot\parallel v\parallel_{Z}\\
  &=\left(a_{k}\ast b_{l'}\right) _{l}\cdot\parallel v\parallel_{Z},
 \end{split}
\end{equation*}
where $a_{k}=2^{k(\sigma+1)}\chi_{\{k;\mid k\mid\leq 4\}}(k)$ and $b_{l'}=2^{l'(\sigma+2)}\parallel\dot{\Delta}_{l'}w\parallel_{L^{1}\left(I;{\cal M}_{\vec{q},\vec{\lambda}}(\mathbb{R}^{d})\right)}$. By applying the Young's inequality for series, we have
\begin{equation*}
\begin{split}
 \parallel  2^{l(\sigma+1)}\parallel I^{l}_{1}\parallel_{L^{1}\left(I;{\cal M}_{\vec{q},\vec{\lambda}}(\mathbb{R}^{d})\right)}\parallel_{l^{r}(\mathbb{Z})}&\leq
 \parallel a_{k}\parallel_{l^{1}(\mathbb{Z})}\parallel b_{l'}\parallel_{l^{r}(\mathbb{Z})}\parallel v\parallel_{Z}\\
 &\leq \parallel w\parallel_{{\cal L}^{1}\left(I;\dot{{\cal N}}^{\sigma+2}_{\vec{q},\vec{\lambda},r}(\mathbb{R}^{d})\right)} \parallel v\parallel_{Z},
\end{split}    
\end{equation*}
which implies the desired estimate for $I^{l}_{1}$. 
A similar procedure is enough to prove that
\begin{equation*}
\parallel  2^{l(\sigma+1)}\parallel I^{l}_{2}\parallel_{L^{1}\left(I;{\cal M}_{\vec{q},\vec{\lambda}}(\mathbb{R}^{d})\right)}\parallel_{l^{q}(\mathbb{Z})}\leq \parallel w\parallel_{Z}\parallel v\parallel_{Z}.    
\end{equation*}
Estimates for the term $I^{l}_{3}$: By the Bernstein's type inequality Lemma 2.4 (i.) (\ref{properties-Besov-spaces-mixed}), we obtain
\begin{equation*}
 \begin{split}
 \parallel I^{l}_{3}\parallel_{L^{1}\left(I;{\cal M}_{\vec{q},\vec{\lambda}}(\mathbb{R}^{d})\right)}&\leq 
 \displaystyle{\sum_{l'\geq l-2}}\parallel \dot{\Delta}_{l'}v\parallel_{L^{1}\left(I;{\cal M}_{\vec{q},\vec{\lambda}}(\mathbb{R}^{d})\right)} \left(
 \displaystyle{\sum_{\mid l''-l'\mid\leq 1}}2^{l''(\sigma+1))}\parallel \dot{\Delta}_{l''}w\parallel_{L^{\infty}\left(I;{\cal M}_{\vec{q},\vec{\lambda}}(\mathbb{R}^{d})\right)}\right)\\
 &=\displaystyle{\sum_{l'\geq l-2}}\parallel \dot{\Delta}_{l'}v\parallel_{L^{1}\left(I;{\cal M}_{\vec{q},\vec{\lambda}}(\mathbb{R}^{d})\right)} \left(
 \displaystyle{\sum_{\mid l''-l'\mid\leq 1}}2^{l''}\cdot 2^{l''\sigma}\parallel \dot{\Delta}_{l''}w\parallel_{L^{\infty}\left(I;{\cal M}_{\vec{q},\vec{\lambda}}(\mathbb{R}^{d})\right)}
 \right).
 \end{split}   
\end{equation*}
By the H\"older's inequality for series we obtain
\begin{equation*}
 \begin{split}
 \parallel I^{l}_{3}\parallel_{L^{1}\left(I;{\cal M}_{\vec{q},\vec{\lambda}}(\mathbb{R}^{d})\right)}&\leq 
 \displaystyle{\sum_{l'\geq l-2}}\parallel \dot{\Delta}_{l'}v\parallel_{L^{1}\left(I;{\cal M}_{\vec{q},\vec{\lambda}}(\mathbb{R}^{d})\right)} \left(\displaystyle{\sum_{\mid l''-l'\mid \leq 1}} 2^{l''r'}\right)^{\frac{1}{r'}}\parallel w\parallel_{{\cal L}^{\infty}\left(I;\dot{{\cal N}}^{\sigma}_{\vec{q},\vec{\lambda},r}\right)}\\
&=\displaystyle{\sum_{l'\geq l-2}}\parallel \dot{\Delta}_{l'}v\parallel_{L^{1}\left(I;{\cal M}_{\vec{q},\vec{\lambda}}(\mathbb{R}^{d})\right)} 2^{l'} \parallel w\parallel_{{\cal L}^{\infty}\left(I;\dot{{\cal N}}^{\sigma}_{\vec{q},\vec{\lambda},r}\right)}\\
&=\displaystyle{\sum_{l'\geq l-2}}2^{-l'(\sigma+1)}\cdot 2^{l'(\sigma+2)}\parallel \dot{\Delta}_{l'}v\parallel_{L^{1}\left(I;{\cal M}_{\vec{q},\vec{\lambda}}(\mathbb{R}^{d})\right)}\parallel w\parallel_{Z}.
 \end{split}   
\end{equation*}
Multiplying by $2^{l(\sigma+1)}$ we obtain
\begin{equation*}
 \begin{split}
  2^{l(\sigma+1)} \parallel I^{l}_{3}\parallel_{L^{1}\left(I;{\cal M}_{\vec{q},\vec{\lambda}}(\mathbb{R}^{d})\right)}&\leq \displaystyle{\sum_{l'}}
  2^{(l-l')(\sigma+1)}\chi_{\{k; k\leq 2\}}(l-l')\cdot 2^{l'(\sigma+2)}\parallel \dot{\Delta}_{l'}v\parallel_{L^{1}\left(I;{\cal M}_{\vec{q},\vec{\lambda}}(\mathbb{R}^{d})\right)}\parallel w\parallel_{Z}\\
  &=\left(a_{k}\ast b_{l'}\right)_{l},
 \end{split}   
\end{equation*}
where $a_{k}=2^{k(\sigma+1)}\chi_{\{k;k\leq 2\}}(k)$ and $b_{l'}=2^{l'(\sigma+2)}\parallel \dot{\Delta}_{l'}v\parallel_{L^{1}\left(I;{\cal M}_{\vec{q},\vec{\lambda}}(\mathbb{R}^{d})\right)}$.
Now we apply the Young's inequality for series to get
\begin{equation*}
\begin{split}
 \parallel   2^{l(\sigma+1)} \parallel I^{l}_{3}\parallel_{L^{1}\left(I;{\cal M}_{\vec{q},\vec{\lambda}}(\mathbb{R}^{d})\right)}\parallel_{l^{r}(\mathbb{Z})}&\leq 
 \parallel a_{k}\parallel_{l^{1}(\mathbb{Z})}\parallel b_{l'}\parallel_{l^{r}(\mathbb{Z})}\parallel w\parallel_{Z}\\
 &\leq \parallel v\parallel_{{\cal L}^{1}\left(I;\dot{{\cal N}}^{\sigma+2}_{\vec{q},\vec{\lambda}.r}\right)}\parallel w\parallel_{Z},
 \end{split}
 \end{equation*}
which implies the desired estimate, since the sequence $(a_{k})\in l^{1}(\mathbb{Z})$ whenever $\sigma+1>0$.
 
\end{proof}
%%%%%%%%%%%%%%%%%%%%%%%%%%%%%%%%%%%%
%%%%%%%%%%%%%%%%%%%%%%%%%%%%%%%%%%
\begin{lema}
 Let us consider $\vec{q}\in[1,\infty)^{d}$ and $\vec{\lambda}\in[0,1)^{d}$ be such that $\displaystyle{\sum_{i=1}^{d}}\left(1-\frac{1-\lambda_{i}}{q_{i}}\right)>0$, $r\in[1,\infty]$ and $s=-1+\displaystyle{\sum_{i=1}^{d}}\left(1-\frac{1-\lambda_{i}}{q_{i}}\right)$. Then, there exists a constant $K_{0}>0$, independent on $\nu$, such that
 \begin{equation*}
  \parallel B(v,w)\parallel_{Z}\leq K_{0}\max{(1,\frac{1}{\nu})}\parallel v\parallel_{Z}\parallel w\parallel_{Z},   
 \end{equation*}
 for each $v,w\in Z={\cal L}^{\infty}\left(I;\dot{{\cal FN}}^{s}_{\vec{q},\vec{\lambda},r}(\mathbb{R}^{d})\right)\bigcap {\cal L}^{1}\left(I;\dot{{\cal FN}}^{s+2}_{\vec{q},\vec{\lambda},r}(\mathbb{R}^{d})\right)$.
\end{lema}
\begin{proof}
%%%%%%%%%%%%%%%%%%
First, we observe that since $B(v,w)(t,\cdot)=A(\mbox{div}(v\otimes w))(t,\cdot)$, where the linear auxiliar operator $A(\cdot)$ is given in (\ref{auxiliar-linearop}), we get 
\begin{equation*}
 \parallel B(v,w)\parallel_{Z}\leq c\max{(1,\frac{1}{\nu})}\parallel \mbox{div}\left(v\otimes w\right)\parallel_{{\cal L}^{1}\left(I;\dot{{\cal FN}}^{s}_{\vec{q},\vec{\lambda},r}(\mathbb{R}^{d})\right)}\leq c\max{(1,\frac{1}{\nu})}\parallel v\otimes w\parallel_{{\cal L}^{1}\left(I;\dot{{\cal FN}}^{s+1}_{\vec{q},\vec{\lambda},r}(\mathbb{R}^{d})\right)},   
\end{equation*}
where $c$ is a constant which does not depend on $\nu$. Therefore, it is enought o prove that 
\begin{equation*}
 \parallel v\otimes w\parallel_{{\cal L}^{1}\left(I;\dot{{\cal FN}}^{s+1}_{\vec{q},\vec{\lambda},r}(\mathbb{R}^{d})\right)}\leq \parallel v\parallel_{Z}\parallel w\parallel_{Z}.
\end{equation*}
%%%%%%%%%%%%%%%%
By the Bony's decomposition we have
\begin{equation*}
\begin{split}
\phi_{l}\left(vw\right)^{\wedge}&=\displaystyle{\sum_{\mid l-l'\mid\leq 4}}\phi_{l}[\left(\dot{S}_{l'-1}v\right)^{\wedge}\ast\left(\phi_{l'}\hat{w}\right)] + \displaystyle{\sum_{\mid l-l'\mid\leq 4}}\phi_{l}[\left(\dot{S}_{l'-1}w\right)^{\wedge}\ast\left(\phi_{l'}\hat{v}\right)] + \displaystyle{\sum_{l'\geq l-2}}\phi_{l}[\left(\phi_{l'}\hat{v}\right)\ast\left( \tilde{\phi}_{l'}\hat{w}\right)]\\
&=I^{l}_{1}+I^{l}_{2}+I^{l}_{3}.
\end{split}
\end{equation*}
Estimates for $I^{l}_{1}$: By the Bernstein's type inequality Lemma 2.6 (i.) (\ref{properties-Fourier-Besov-mixed}), $\parallel\phi_{l''}\hat{v}\parallel_{L^{1}(\mathbb{R}^{d})}\leq 2^{l''(s+1)}\parallel\phi_{l''}\hat{v}\parallel_{{\cal M}_{\vec{q},\vec{\lambda}}(\mathbb{R}^{d})}$, we have
\begin{equation*}.
\begin{split}
\parallel I^{l}_{1}\parallel_{L^{1}\left(I;{\cal M}_{\vec{q},\vec{\lambda}}(\mathbb{R}^{d})\right)}&\leq \displaystyle{\sum_{\mid l-l'\mid\leq 4}}\parallel \parallel \left(S_{l'-1}v\right)^{\wedge}\ast \left(\phi_{l'}\hat{w}\right)\parallel_{{\cal M}_{\vec{q},\vec{\lambda}}(\mathbb{R}^{d})}\parallel_{L^{1}(I)}\\
& \leq \displaystyle{\sum_{\mid l-l'\mid\leq 4}}\parallel \parallel \left(S_{l'-1}v\right)^{\wedge}\parallel_{L^{1}(\mathbb{R}^{d})} \parallel \phi_{l'}\hat{w}\parallel_{{\cal M}_{\vec{q},\vec{\lambda}}(\mathbb{R}^{d})}\parallel_{L^{1}(I)}\\
&\leq\displaystyle{\sum_{\mid l-l'\mid\leq 4}}
\left(\displaystyle{\sum_{l''\leq l'}}2^{l''(s+1)} \parallel \phi_{l''}\hat{v}\parallel_{L^{\infty}(I;{\cal M}_{\vec{q},\vec{\lambda}}(\mathbb{R}^{d}))}\right) \parallel \phi_{l'}\hat{w}\parallel_{L^{1}(I;{\cal M}_{\vec{q},\vec{\lambda}}(\mathbb{R}^{d}))}.
\end{split}    
\end{equation*}
Then,
\begin{equation*}
\parallel I^{l}_{1}\parallel_{L^{1}\left(I;{\cal M}_{\vec{q},\vec{\lambda}}(\mathbb{R}^{d})\right)}\leq 
\displaystyle{\sum_{\mid l-l'\mid\leq 4}}
\left( \displaystyle{\sum_{l''\leq l'}}2^{l''}2^{sl''}
\parallel \phi_{l''}\hat{v}\parallel_{L^{\infty}\left((I;{\cal M}_{\vec{q},\vec{\lambda}}(\mathbb{R}^{d})\right)}\right) \parallel \phi_{l'}\hat{w}\parallel_{L^{1}\left(I;{\cal M}_{\vec{q},\vec{\lambda}}(\mathbb{R}^{d})\right)}.
\end{equation*}
\end{proof}
By H\"older's inequality for series we get
\begin{equation*}
\parallel I^{l}_{1}\parallel_{L^{1}\left(I;{\cal M}_{\vec{q},\vec{\lambda}}(\mathbb{R}^{d})\right)}\leq \displaystyle{\sum_{\mid l-l'\mid\leq 4}}\left(\displaystyle{\sum_{l''\leq l'}}2^{l''r'} \right)^{\frac{1}{r'}}\cdot\parallel \phi_{l'}\hat{w}\parallel_{L^{1}\left(I;{\cal M}_{\vec{q},\vec{\lambda}}(\mathbb{R}^{d})\right)} \parallel v\parallel_{{\cal L}^{\infty}\left(I;\dot{{\cal FN}}^{s}_{\vec{q},\vec{\lambda},r}\left(\mathbb{R}^{d}\right)\right)}.
\end{equation*}
If we multiply by $2^{l(s+1)}$ in both sides we have
\begin{equation*}
\begin{split}
2^{l(s+1)}\parallel I^{l}_{1}\parallel_{L^{1}\left(I;{\cal M}_{\vec{q},\vec{\lambda}}(\mathbb{R}^{d})\right)}&\leq
\left(
\displaystyle{\sum_{l'}}2^{(l-l')(s+1)}\chi_{\{k;\mid k\mid\leq 4\}}\left(l-l'\right)2^{l'(s+2)}\parallel \phi_{l'}\hat{w}\parallel_{L^{1}\left(I;{\cal M}_{\vec{q},\vec{\lambda}}(\mathbb{R}^{d})\right)}
\right)\cdot \parallel v\parallel_{Z}\\
&=\left(a_{k}\ast b_{l'}\right)_{l}\cdot \parallel v\parallel_{Z},
\end{split}    
\end{equation*}
where $a_{k}=2^{k(s+1)}\chi_{\{k;\mid k\mid\leq 4\}}\left(k\right)$ and $b_{l'}=2^{l'\left(s+2\right)}\parallel \phi_{l'}\hat{w}\parallel_{L^{1}\left(I;{\cal M}_{\vec{q},\vec{\lambda}}(\mathbb{R}^{d})\right)}$. By Young's inequality for series, we obtain
\begin{equation*}
\begin{split}
\parallel 2^{l(s+1)}\parallel I^{l}_{1}\parallel_{L^{1}\left(I;{\cal M}_{\vec{q},\vec{\lambda}}(\mathbb{R}^{d}) \right)}\parallel_{l^{r}\left(\mathbb{Z}\right)}&\leq   \parallel a_{k}\parallel_{l^{1}\left(\mathbb{Z}\right)}\cdot 
\parallel b_{l'}\parallel_{l^{r}\left(\mathbb{Z}\right)}\cdot \parallel v\parallel_{Z}\\
&\leq \parallel w\parallel_{{\cal L}^{1}\left(I;\dot{{\cal FN}}^{s}_{\vec{q},\vec{\lambda},r}\right)}\cdot \parallel v\parallel_{Z}\\
&\leq \parallel w\parallel_{Z}\cdot \parallel v\parallel_{Z}.
\end{split}    
\end{equation*}
Estimates for $I^{l}_{2}$: Similar computations as for $I^{l}_{1}$ provides
\begin{equation*}
\parallel 2^{l(s+1)}\parallel I^{l}_{2}\parallel_{L^{1}\left(I;{\cal M}_{\vec{q},\vec{\lambda}}(\mathbb{R}^{d}) \right)}\parallel_{l^{r}\left(\mathbb{Z} \right)}\leq \parallel w\parallel_{Z} \parallel v\parallel_{Z}.     
\end{equation*}
Estimates for the term $I^{l}_{3}$: Since $\parallel \phi_{l}\hat{u}\parallel_{L^{1}(\mathbb{R}^{d})}\leq 2^{l(s+1)}\parallel \phi_{l}\hat{u}\parallel_{{\cal M}_{\vec{q},\vec{\lambda}}(\mathbb{R}^{d})}$, we get
\begin{equation*}
\begin{split}
 \parallel I^{l}_{3}\parallel_{L^{1}\left(I;{\cal M}_{\vec{q},\vec{\lambda}}(\mathbb{R}^{d})\right)}&\leq 
 \displaystyle{\sum_{l'\geq l-2}}\left(
 \displaystyle{\sum_{\mid l''-l'\mid\leq 1}}2^{l''(s+1))}\parallel \phi_{l''}\hat{w}\parallel_{L^{\infty}\left(I;{\cal M}_{\vec{q},\vec{\lambda}}(\mathbb{R}^{d})\right)}
 \right) 
 \parallel \phi_{l'}\hat{v}\parallel_{L^{1}\left(I;{\cal M}_{\vec{q},\vec{\lambda}}(\mathbb{R}^{d})\right)}\\
 &\leq  \displaystyle{\sum_{l'\geq l-2}}\left(
 \displaystyle{\sum_{\mid l''-l'\mid\leq 1}}2^{l''}\cdot 2^{s l''}\parallel \phi_{l''}\hat{w}\parallel_{L^{\infty}\left(I;{\cal M}_{\vec{q},\vec{\lambda}}(\mathbb{R}^{d})\right)}
 \right) 
 \parallel \phi_{l'}\hat{v}\parallel_{L^{1}\left(I;{\cal M}_{\vec{q},\vec{\lambda}}(\mathbb{R}^{d})\right)}.
\end{split}    
\end{equation*}
By the H\"older's inequality for series, we have
\begin{equation*}
\begin{split}
  \parallel I^{l}_{3}\parallel_{L^{1}\left(I;{\cal M}_{\vec{q},\vec{\lambda}}(\mathbb{R}^{d})\right)}&\leq
  \displaystyle{\sum_{l'\geq l-2}}
  \left(\sum_{\mid l''-l'\mid\leq 1}2^{l''r'}\right)^{\frac{1}{r'}} \parallel \phi_{l'}\hat{v}\parallel_{L^{1}\left(I;{\cal M}_{\vec{q},\vec{\lambda}}(\mathbb{R}^{d})\right)}\cdot\parallel w\parallel_{{\cal L}^{\infty}\left(I;\dot{{\cal FN}}^{s}_{\vec{q},\vec{\lambda},r}(\mathbb{R}^{d})\right)}\\
  &\leq   \displaystyle{\sum_{l'\geq l-2}} 2^{-l'(s+1)}\left( 2^{l'(s+2)}\parallel \phi_{l'}\hat{v}\parallel_{L^{1}\left(I;{\cal M}_{\vec{q},\vec{\lambda}}(\mathbb{R}^{d})\right)}\right) \parallel w\parallel_{Z}
\end{split}    
\end{equation*}
By multiplying both sides by $2^{l(s+1)}$, we obtain
\begin{equation*}
\begin{split}
 2^{l(s+1)}  \parallel I^{l}_{3}\parallel_{L^{1}\left(I;{\cal M}_{\vec{q},\vec{\lambda}}(\mathbb{R}^{d})\right)}&\leq \displaystyle{\sum_{l'}}
 2^{\left(l-l'\right)\left(s+1\right)}\chi_{\{k;k\leq 2\}]}(l-l')\cdot \left( 2^{l'(s+2)}\parallel \phi_{l'}\hat{v}\parallel_{L^{1}\left(I;{\cal M}_{\vec{q},\vec{\lambda}}(\mathbb{R}^{d})\right)}\right)\cdot \parallel w\parallel_{Z}\\
 &\leq \left(a_{k}\ast b_{l'}\right)_{l}\cdot \parallel w\parallel_{Z},
\end{split}    
\end{equation*}
where $a_{k}=2^{k(s+1)}\chi_{\{k;k\leq 2\}}(k)$ and $b_{l'}=2^{l'(s+2)}\parallel \phi_{l'}\hat{v}\parallel_{L^{1}\left(I;{\cal M}_{\vec{q},\vec{\lambda}}(\mathbb{R}^{d})\right)}$.
Applying the Young's inequality for series, we obtain
\begin{equation*}
\begin{split}
\parallel 2^{l(s+1)}  \parallel I^{l}_{3}\parallel_{L^{1}\left(I;{\cal M}_{\vec{q},\vec{\lambda}}(\mathbb{R}^{d})\right)} \parallel_{l^{r}\left(\mathbb{Z}\right)}&\leq 
\parallel a_{k}\parallel_{l^{1}(\mathbb{Z})}\cdot\parallel b_{l'}\parallel_{l^{r}(\mathbb{Z})}\cdot \parallel w\parallel_{Z}\\
&\leq \parallel v\parallel_{Z}\cdot \parallel w\parallel_{Z},
\end{split}    
\end{equation*}
where the sequence $(a_{k})_{k\in\mathbb{Z}}\in l^{1}(\mathbb{Z})$ if $s+1>0$. These computations are enough to conclude the proof.

%%%%%%%%%%%%%%%%%%%%%
\subsection{Proof of the last two main theorems}
\label{lastsection}
%%%%%%%%%%%%%%%%%%%%%

In this subsection we conclude the proof of the main theorems. We set the computations for the Besov mixed-Morrey space, since the procedure for Fourier-Besov mixed-Morrey spaces are similar. Let us have in mind expressions (\ref{mild-formulation}), (\ref{mild-formulation-abstract}) and Lemma \ref{fixed-point-scheme}. We first observe that, by Lemma \ref{estimates-lin-besov} we have
\begin{equation*}
 \parallel z_{0}\parallel_{Z}\leq c\max{(1,\frac{1}{\nu})}\parallel u_{0}\parallel_{\dot{{\cal N}}^{\sigma}_{\vec{q},\vec{\lambda},r}(\mathbb{R}^{d})},   
\end{equation*}
where $c> 0$ does not depend on $\nu$. By Lemma \ref{estimates-bilinear-besov} we have
\begin{equation*}
 \parallel B\parallel_{{\cal B}(Z)}\leq K_{0}\max{(1,\frac{1}{\nu})} =:K.   
\end{equation*}
First we take 
\begin{equation*}
0<\varepsilon<\frac{1}{4K}=\frac{1}{4K_{0}\max{(1,\frac{1}{\nu})}}=\frac{\min{(1,\nu)}}{4K_{0}}.    
\end{equation*}
Therefore, in order to get the global well-posedness result it is enough to consider $u_{0}\in \dot{{\cal N}}^{\sigma}_{\vec{q},\vec{\lambda},r}(\mathbb{R}^{d})$ such that
\begin{equation*}
 \parallel u_{0}\parallel_{\dot{{\cal N}}^{\sigma}_{\vec{q},\vec{\lambda},r}(\mathbb{R}^{d})}\leq\frac{\varepsilon}{c\max{(1,\frac{1}{\nu})}}=\frac{\varepsilon}{c}\min{(1,\nu)}.   
\end{equation*}
This conclude the proof of item (i.) from the first main theorem. For the proof of the second theorem we proceed similarly. Finally, the procedure for the proof of the weak continuity from $[0,\infty)$ toward ${\cal S}'(\mathbb{R}^{d})$ is well-known and we omit it here.

%%%%%%%%%%%%%%%%%%%%
\section*{Data availability statement}
%%%%%%%%%%%%%%%%%%%%
Data sharing is not applicable to this article as no data sets were 
generated or analysed during the current study.

%%%%%%%%%%%%%%%%%%%%%%%%%%%%%%
 \section*{Conflict of Interest}
%%%%%%%%%%%%%%%%%%%%%%%%
There is no conflict of interests.

\end{document}